\documentclass[12pt]{amsart}

\usepackage{graphicx,psfrag,epsfig}
\usepackage{amsmath}
\usepackage{amsthm}
\usepackage{bbm}
\usepackage[left=1.25in, right=1.25in, top=1.25in, bottom=1.25in]{geometry}
\usepackage{amsfonts}
\usepackage{amssymb}
\usepackage{verbatim}
\usepackage{comment}
\usepackage[pdftex]{xcolor}
\usepackage{mathtools}
\usepackage{ulem}
\usepackage[parfill]{parskip}
\usepackage{cleveref}

\newtheorem{theorem}{Theorem}
\newtheorem{lemma}{Lemma}[section]
\newtheorem{sublemma}[lemma]{Sublemma}
\newtheorem{corollary}[lemma]{Corollary}
\newtheorem{proposition}[lemma]{Proposition}

\theoremstyle{definition}

\newtheorem{remark}[lemma]{Remark}

\numberwithin{equation}{section}

\newcommand{\ol}[1]{\overline{#1}}

\everymath{\displaystyle}

\def\fp{{\mathfrak{p}}}

\def\hcX{{\mathcal{X}}}
\def\brcL{{\mathcal{L}}}

\def\bareps{{\bar{\eps}}}

\def\DS{\displaystyle}

\def\Phase{{\mathbb{M}}}

\newcommand{\KF}[1]{{\textcolor{red}{#1}}}

\def\dist{{\rm dist}}
\def\eps{{\varepsilon}}

\def\Card{{\rm Card}}

\def\Cov{{\rm Cov}}

\def\Prob{{\mathbb{P}}}

\def\Vol{{\rm Vol}}
\def\EXP{{\mathbb{E}}}

\def\naturals{\mathbb{N}}

\def\Tor{\mathbb{T}}

\def\reals{\mathbb{R}}

\def\integers{\mathbb{Z}}

\def\bE{\mathbf{E}}

\def\bP{\mathbf{P}}

\def\ba{\mathbf{a}}

\def\bb{\mathbf{b}}

\def\be{\mathbf{e}}

\def\bk{\mathbf{k}}

\def\bm{\mathbf{m}}

\def\bp{\mathbf{p}}
\def\bq{\mathbf{q}}

\def\bu{\mathbf{u}}

\def\bw{\mathbf{w}}
\def\bx{\mathbf{x}}

\def\bnu{\boldsymbol{\nu}}

\def\btheta{{\boldsymbol{\theta}}}
\def\bomega{{\boldsymbol{\omega}}}

\def\bDelta{{\boldsymbol{\Delta}}}

\def\brC{{\bar C}}

\def\brcX{{\bar\cX}}

\def\brc{{\bar c}}

\def\brm{{\bar m}}

\def\brs{{\bar s}}

\def\brDelta{{\bar \Delta}}

\def\brPhi{{\bar \Phi}}

\def\cB{\mathcal{B}}

\def\cI{\mathcal{I}}

\def\cE{\mathcal{E}}

\def\cL{\mathcal{L}}

\def\cM{\mathcal{M}}

\def\cO{\mathcal{O}}

\def\cS{\mathcal{S}}

\def\cU{\mathcal{U}}

\def\cV{\mathcal{V}}

\def\cX{\mathcal{X}}

\def\cY{\mathcal{Y}}
\def\cZ{\mathcal{Z}}

\def\fA{\mathfrak{A}}

\def\fB{\mathfrak{B}}

\def\fN{\mathfrak{N}}

\def\fa{\mathfrak{a}}

\def\fn{\mathfrak{n}}

\def\hcE{{\hat{\cE}}}

\def\hchi{{\hat\chi}}

\def\hDelta{\hat{\mathbf{\Delta}}}

\def\tcI{{\tilde\cI}}

\def\tcX{{\tilde \cX}}

\def\tc{{\tilde c}}

\def\tp{{\tilde p}}

\def\tDelta{{\tilde{\mathbf{\Delta}}}}

\def\tchi{{\tilde\chi}}

\def\tGamma{{\tilde{\Gamma}}}

\makeatother

\def\beq{\begin{equation}}
\def\eeq{\end{equation}}

\subjclass[2000]{Primary 60F05; Secondary 11J13, 37A50}
\begin{document}

\author[Dmitry Dolgopyat]{Dmitry Dolgopyat$^{1*}$}
\author[Kasun Fernando]{Kasun Fernando$^2$}
\title[An error term in the CLT for discrete random variables.]
      {An error term in the Central Limit Theorem for sums of discrete random
        variables.}

\maketitle
(*) corresponding author

(1) Department of Mathematics, University of Maryland, College Park MD, USA
{\em email:} dmitry@math.umd.edu

(2) Centro De Giorgi, Scuola Normale Superiore, Pisa PI 56126, Italy \\
{\em email:}  buddhima.akurugodage@sns.it

\begin{abstract}
We consider sums of independent identically distributed random variables whose distributions have $d+1$ atoms. Such distributions never admit an Edgeworth expansion of order $d$ but we show that for almost all parameters the Edgeworth expansion of order $d-1$ is valid and the error of the order $d-1$ Edgeworth expansion is typically of order
$n^{-d/2}.$
\end{abstract}

\section{Introduction.}
Let $X$ be a random variable with zero mean and variance $\sigma^2.$
Let $\DS S_n=\sum_{j=1}^n X_j$ where $X_j$ are independent
identically distributed and have the same distribution as $X.$ The Central Limit Theorem says that for each~$z$
$$ \lim_{n\to\infty} \Prob\left(\frac{S_n}{\sigma \sqrt{n}}\leq z\right)=\fN(z) $$
where $$\fN(z)=\int_{-\infty}^z \fn (y) dy \text{ and } \fn(y)=\frac{1}{\sqrt{2\pi} } e^{-y^2/2}. $$
A classical problem in probability theory is computing higher order approximations to
$\DS \Prob\left(\frac{S_n}{\sigma \sqrt{n}}\leq z\right)$.
In particular, the
 order $r$ Edgeworth series of $S_n$ is an expression of the form
$$ \cE_r (z)=\fN(z)+\fn(z) \sum_{k=1}^r \frac{P_k(z)}{n^{k/2}} $$
where $P_k$ are polynomials
such that the characteristic function $\phi(t)=\EXP(e^{it X})$ and the Fourier transform $\hcE_r$ of $\cE_r$
satisfy
$$ \phi\left(\frac{t}{\sigma \sqrt{n}}\right)^n-\hcE_r(t)=o\left(n^{-r/2}\right). $$
In particular,
$$ \cE_1(z)=\fN(z)+\fn(z) \frac{\EXP(X^3)}{6\sigma^3\sqrt{n}} (1-z^2), $$
\begin{multline*}
\cE_2(z)=\fN(z)+\fn(z) \left[\frac{\EXP(X^3)}{6\sqrt{n}\sigma^3} (1-z^2)+\frac{\EXP(X^4)-3\sigma^4}{24 n \sigma^4} (3z-z^3) \right. \\ \left. -
\frac{\EXP(X^3)^2}{72 n\sigma^6}
(15z-10 z^3+z^5)
\right].
\end{multline*}
We say that $S_n$ admits an order $r$ Edgeworth expansion if for all $z$
\begin{equation}
\label{EdgeValid}
 \lim_{n\to \infty} n^{r/2} \left[\Prob\left(\frac{S_n}{\sigma \sqrt{n}}\leq z\right)-\cE_r(z)\right]=0.
\end{equation}

{Recall that a lattice random variable is a discrete random variable taking values on a set of points of the form $a+nh, n \in \integers$, where $h>0$, $a\in \reals$}. It is known that $S_n$ admits the first order Edgeworth expansion if and only if $X$ is non-lattice (see \cite{Es}). The problem of higher order expansion is more complicated. For example, a sufficient condition for $S_n$ to admit the order $r$ Edgeworth expansion is that $\EXP(|X|^{r+2})<\infty$  and $X$ has a density. But this condition is far from necessary. We refer the reader to \cite[Chapter XVI]{Fel} for discussion of these and
related results. We also note that \cite{Ba, B} discusses a weak Edgeworth expansion where the LHS of \eqref{EdgeValid} is convolved with smooth compactly supported functions.

In this paper, we consider a case which is opposite to $X$ having a density, namely we suppose that $X$ has a discrete distribution with $d+1$ atoms where $d\geq 2.$
$d=2$ is the simplest non-trivial case since the distributions with two atoms are lattice, and as a result, they do not admit even the first order Edgeworth expansion.

Thus we suppose that $X$ takes values $a_1, \dots, a_{d+1}$ with probabilities $p_1, \dots, p_{d+1}$, respectively. Since $X$ should have zero mean we suppose that our $2(d+1)-$tuple $(\ba, \bp)$ belongs to the set
$$ \Omega=\{p_i>0,\quad p_1+\dots+p_{d+1}=1, \quad p_1 a_1+\dots+p_{d+1} a_{d+1}=0\}  .$$
It is easy to see that $S_n$ never admits the order $d$ Edgeworth expansion.
Indeed,
\begin{equation}
\label{PdfSum}
\Prob_{\ba, \bp}(S_n \leq z)=
\sum_{{m_i \geq 0, \;\sum m_i=n}\atop{\sum m_ia_i\leq z}}
\frac{n!}{m_1! \dots m_{d+1}!} p_1^{m_1} \dots  p_{d+1}^{m_{d+1}} .
\end{equation}
The Local Central Limit Theorem (see \cite[Theorem 2.1.1]{Law}), applied to the time homogeneous $\integers^d$-random walk which jumps to $\be_i$ from the origin $\boldsymbol{0}$ with probability $p_i$ for $i=1,\dots, d$ and stays at $\boldsymbol{0}$ with probability $p_{d+1},$
gives us that for all $\eps$ there is $n_0\in\naturals$ such that for all $n\geq n_0$
and all $\bm \in \integers^d$
\begin{equation}\label{EstLCLT}
  {\left|\Prob(T_n=\bm)-\frac{1}{\sqrt{(2\pi n)^d \, \text{det} \Gamma } }e^{-\frac{(\bm-n\bq) . \Gamma^{-1} (\bm-n\bq)}{2n}} \right| \leq  \frac{\eps}{n^{d/2}}}
\end{equation}
where $T_n$ is the position of the random walk after $n$ steps, $\Gamma$ is the associated covariance matrix and $\bq=(p_1,\dots,p_{d})$.
Also, if $m_1,\dots m_d, m_{d+1}$ are integers such that
$m_1+\dots+m_d+m_{d+1}=n$ and $m_i \geq 0$, then, taking $\bm=(m_1,\dots,m_d),$ we have
$$\Prob(T_n=\bm)=\frac{n!}{m_1! \dots m_{d+1}!}p^{m_1}_1 \dots p^{m_{d+1}}_{d+1}.$$

As a result, if
$$\sum m_ia_i = n\sum a_ip_i +\mathcal{O}(\sqrt{n}),$$
then {in \eqref{EstLCLT}, the exponent of $e$ is $\cO(1)$, and hence, for sufficiently large $n$,}
$$n^{d/2}\Prob(T_n=\bm) = n^{d/2}\frac{n!}{m_1! \dots m_{d+1}!}p^{m_1}_1 \dots p^{m_{d+1}}_{d+1}$$
is uniformly bounded from below. 
Accordingly, {from \eqref{PdfSum}}, it follows that $\Prob_{\ba, \bp}(S_n \leq z)$ has jumps of order $n^{-d/2}.$ On the other hand, $\cE_d(z)$ is a smooth function of $z$. So, it can not approximate both $\Prob_{\ba, \bp}(S_n \leq z-0)$ and $\Prob_{\ba, \bp}(S_n \leq z+0)$ at the points of jumps {without making an error of $\cO(n^{-d/2})$. This means that it is {\it not} true that $\Prob_{\ba, \bp}(S_n \leq z)=\cE_d(z)+o(n^{-d/2})$ for all $z$, showing that the order $d$ Edgeworth expansion fails.}

{However, in this paper, we show that for \textit{typical} $(\ba, \bp)$ this failure of the order $d$ Edgeworth expansion happens just barely}. We present two results in this direction. For the first result, let
$$ b_j=a_j-a_1, \text{ for }j=2,\dots, d+1. $$
Then, the characteristic function of $X$, $\phi$, satisfies
$$ \phi(s)=e^{is a_1}\psi(s)\quad\text{where}\quad
\psi(s)=p_1+p_2 e^{is b_2}+\dots+p_{d+1} e^{is b_{d+1}}.$$
Set
$$d(s)=\max_{j\in \{2,\dots, d+1\}} \dist(b_j s, 2\pi \integers). $$
We say that $\ba$ is $\beta$-Diophantine if there is a constant $K$ such that for $|s|>1$,
$$ d(s)\geq \frac{K}{|s|^\beta}. $$
{It follows from the classical Khinchine-Groshev Theorem  (see e.g. \cite[Theorem 1.1]{KM2}
or \cite{Sp})} that almost every $\ba$ is $\beta$-Diophantine provided that $\beta>\frac{1}{d-1}.$
\begin{theorem}
  \label{ThDioEdge}
  If $\ba$ is $\beta$-Diophantine and
  \begin{equation}
   \label{RSmall}
   2\left(R-\frac{1}{2}\right)\beta<1
   \end{equation}
  then
$$ \lim_{n\to\infty} n^R \left[\Prob_{\ba, \bp} \left(\frac{S_n}{\sigma \sqrt{n}}\leq z\right)-\cE_{d-1}(z)\right]
=0$$
{uniformly in $z \in \reals$.}
\end{theorem}
\noindent
Thus, for almost every $\ba$ the order $(d-1)$ Edgeworth expansion approximates the distribution of
$\frac{S_n}{\sigma \sqrt{n}}$ with error $\cO(n^{\eps-d/2})$ for any $\varepsilon.$

Note that \Cref{ThDioEdge} applies for all $\beta$s, and in particular, for $\beta$s which are much larger than $\frac{1}{d-1}.$ However, if $\beta$ is large, then the statement of the theorem can be simplified. Namely, let $r$ be the integer such that $r<2R\leq r+1.$ (Note that \eqref{RSmall} can be rewritten as
$2R<\frac{1}{\beta}+1$. So, provided that $2R$ is sufficiently close to $\frac{1}{\beta}+1$ we can take
$r=\left\langle \frac{1}{\beta}\right\rangle+1$ where $\langle s \rangle$ denotes the largest integer which is strictly smaller than $s$.) Then,
\begin{align*}
\Prob_{\ba, \bp}\left(\frac{S_n}{\sigma \sqrt{n}}\leq z\right) &=\cE_{d-1}(z)+{ o\left(\frac{1}{n^R}\right)} \\ &= \cE_{r}(z)+{  o\left(\frac{1}{n^R}\right)}+\cO\left(\cE_{d-1}(z)-\cE_{r}(z) \right).
\end{align*}
Since $\frac{r+1}{2}>R$ the first error term dominates the second and we obtain the following result.
\begin{corollary} \label{ThDioEdgeCor} {Suppose that $\ba$ is $\beta$-Diophantine, $r=1+\left\langle\frac{1}{\beta}\right\rangle,$ and $r<2R<\frac{1}{\beta}+1.$ Then}
$$ \lim_{n\to\infty} n^R \left[\Prob_{\ba, \bp} \left(\frac{S_n}{\sigma \sqrt{n}}\leq z\right)-\cE_{r}(z)\right]
=0 $$
{uniformly in $z \in \reals$}.  
\end{corollary}
\noindent
\Cref{ThDioEdgeCor} shows that
for almost every $\ba$ and for $r \in \{1,\dots,d-1\}$, the order $r$ Edgeworth expansion is valid.
Our next results show that 
\begin{equation}
\label{EdgeErrorD}
\Prob_{\ba, \bp} \left(\frac{S_n}{\sigma \sqrt{n}}\leq z\right)-\cE_d(z)
\end{equation}
is typically of order $\cO(n^{-d/2})$ but the $\cO(n^{-d/2})$ term has wild oscillations. To formulate this result precisely, we suppose that our $2(d+1)$-tuple is chosen at random according to an absolutely continuous distribution $\bP$ on $\Omega.$ Thus, \eqref{EdgeErrorD} becomes a random variable.
\begin{theorem}
\label{Th2EdgeRand}
There exists a smooth function $\Lambda(\ba, \bp)$ such that
for each $z$ the random variable
$$ e^{z^2/2} \frac{n^{d/2}}{\Lambda(\ba, \bp)}
\left[\cE_d(z)-\Prob_{\ba, \bp} \left(\frac{S_n}{\sigma \sqrt{n}}\leq z\right)\right]$$
converges in law to a non-trivial random variable $\cX$ {$($defined below in Lemma \ref{LmThetaLike}$)$
whose distribution is independent of $z$ and $\bP$}.
\end{theorem}

{The formulas for the normalizing factor $\Lambda(\ba, \bp)$ and the limiting random variable
$\cX$ are quite complicated and the next few pages are devoted to their definitions.}

{The normalization is defined as follows:}
\begin{equation}
\label{DefLambda}
  \Lambda(\ba, \bp)=\frac{|a_{d+1}-a_1|}{(2\pi)^{d+\frac{1}{2}}\sqrt{\det(D_{\ba, \bp})}\;\sigma(\ba,\bp) }
\end{equation}
where $\sigma(\ba, \bp)$ denotes the standard deviation of the distribution of the random variable taking value
$a_j$ with probability $p_j$ and
$D_{\ba,\bp}$ is a $(d-1)\times (d-1)$ matrix defined as follows.

\subsection*{{The matrix $D_{\ba,\bp}$}} Fix $p_1, \dots, p_{d+1}$ and consider a map
$$ \zeta(\textbf{y})=\left|\sum_{j=1}^{d+1} p_j e^{i y_j} \right|^2 $$
where $\textbf{y}=(y_1, \dots, y_{d+1}) \in \reals^{d+1}$. Let
$Y$ be a random variable taking values $y_j$ with probability $p_j$. Then, for small $\textbf{y}$
we have
$$ \bE\left(e^{iY}\right)=1-\frac{\bE(Y^2)}{2}+i \bE(Y)+\cO\left(\|\textbf{y}\|^3\right) .$$
Hence,
\begin{equation}\label{charNorm}
\zeta(\textbf{y})=1-\bE(Y^2)+\bE(Y)^2+O\left(\|\textbf{y}\|^3\right)=1-V(Y)+\cO\left(\|\textbf{y}\|^3\right)
\end{equation}
where $V(\cdot)$ is the variance.

Next, consider the quadratic form given by $Q(\textbf{y},\textbf{y})=V(Y(\textbf{y})).$
{Let $\bx, \ba \in \reals^{d+1}$ be fixed}. In order to maximize $s \mapsto \zeta(\bx + s\ba)$, we want to minimize $s\mapsto Q(\bx+s\ba,\bx+s\ba)$.
We have
$$ Q(\bx+s\ba, \bx+s\ba)=Q(\bx,\bx)+2s Q(\bx,\ba)+s^2 Q(\ba,\ba).$$
It follows that the minimum is achieved at $s^*=-\frac{Q(\ba,\bx)}{Q(\ba,\ba)}$ and its value is
$$ D(\bx,\bx)=Q(\bx,\bx)-\frac{Q(\bx,\ba)^2}{Q(\ba,\ba)}=\frac{Q(\bx,\bx) Q(\ba,\ba)-Q(\ba, \bx)^2}{Q(\ba,\ba)}$$
$$=
\frac{V(Y(\bx)) V(X)-\Cov^2(X,Y(\bx))}{V(X)} $$
where $X$ is the random variable taking values $a_j$ with probability $p_j.$
Note that
$D(\bx,\bx)>0$ on the subspace $x_1=x_{d+1}=0$ since
$\Cov^2(X, Y(\bx))=V(X) V(Y(\bx))$ iff
{$Y(\bx)=c_1X+c_2$. Note that the RHS takes $(d+1)$ different values if $c_1\neq 0$
and it takes a single value if $c_1=0.$ On the other hand, the LHS takes at most $d$
different values on $\{x_1=x_{d+1}=0\}$ and it takes a single value only at $\bf 0$. This implies that
$Q(\bx)\neq 0$ unless $\bx=\bf 0\,,$ and hence, $Q$ is non degenerate.}
Then $D_{\ba,\bp}$ is the $(d-1)\times(d-1)$ positive definite matrix such that $\KF{-}4D_{\ba, \bp}$ is the Hessian of $ \reals^{d-1} \ni \tilde{\bx} \mapsto \zeta(0, \tilde\bx, 0).$
The formula for $D_{\ba, \bp}$ will be proven in Section \ref{SS53} (see \eqref{DefD1}).

{We note that the infinitesimal computation described above is relevant because we will show, in the course of proving Theorem \ref{Th2EdgeRand}, that the main contribution to the error term come from the resonant points where the Taylor expansion could be used. See Section~\ref{SS53} for more details.}

To define $\cX$, we need some notation. Let $\Phase$ be the space of pairs $(\cL, \chi)$ where $\mathcal{L}$ is a unimodular lattice in $\reals^d$ and $\chi$ is a character, that is, a
homomorphism $\chi:\cL\to\Tor.$

\subsection*{{The Haar measure on $\Phase$}} The Haar measure $\mu$ on $\Phase$ can be defined in two equivalent ways.
First, note that $\chi$ is of the form $\chi(\bw)=e^{2\pi i \tchi(\bw)}$ for some linear functional
$\tchi\in(\reals^d)^*.$
$SL_d(\reals)$
acts on $\reals^d\oplus (\reals^d)^*$ by the formula
$$ A(\bw, \tchi)=(A\bw, \tchi A^{-1} ). $$
Observe that if $A(\bw, \tchi)=(\hat{\bw}, \hchi)$ then
\begin{equation}
\label{InvProd}
\tchi(\bw)=\hchi(\hat{\bw}).
\end{equation}
The above action of $SL_d(\reals)$ induces the following action of
$SL_d(\reals)\ltimes (\reals^d)^*$ on $\Phase$
$$ (A, \tchi)(\cL, \chi)=(A\cL, e^{2\pi i\tchi} \cdot (\chi\circ A^{-1})). $$
This action is transitive because $SL_d(\reals)$ acts transitively on unimodular lattices and
$(\reals^d)^*$ acts transitively on characters. This allows us to identify $\Phase$ with
$$(SL_d(\reals)\ltimes \reals^d)/(SL_d(\integers)\ltimes \integers^d)$$ and so $\Phase$ inherits the Haar measure from
$SL_d(\reals)\ltimes \reals^d.$

The second way to define the Haar measure
is to note that the space $\cM$ of unimodular lattices is naturally identified with
$SL_d(\reals)/SL_d(\integers)$, and so, it inherits the Haar measure from $SL_d(\reals).$ Next, for a fixed
$\cL$ the set of homomorphisms 
$\chi:\cL\to\Tor$ is a $d$ dimensional torus. So, it comes with its own Haar measure.
Now, if we want to compute the average of a function $\Phi(\cL, \chi)$ with respect to the Haar measure
then we can first compute its average $\brPhi(\cL)$ in each fiber and then integrate the result with respect to
the Haar measure on the space of lattices. In the proof of Lemma \ref{LmThetaLike}
given in Section \ref{ScASConv}, the averaging inside a fiber
will be denoted by $\bE_\chi$ and the averaging with respect to the Haar measure on the space of lattices will be denoted by
$\bE_\cL.$

\subsection*{{The random variable $\cX$}} Given a vector $\bw\in \reals^d$, we denote by $y(\bw)$ its first coordinate and by $\bx(\bw)$ its last $d-1$ coordinates. We also denote by $\|\cdot\|$ the standard Euclidean norm.

\begin{lemma}
\label{LmThetaLike}
For almost every pair $(\cL, \chi)\in \Phase$ with respect to the Haar measure
the following limit exists
\begin{equation}
\label{Lim2k}
\cX(\cL,\chi)=\lim_{R\to\infty}
\sum_{\bw\in\cL \setminus{\{\boldsymbol{0}\},\;\; \|\bw\|\leq R}}
\frac{\sin(2\pi\chi(\bw))}{y(\bw)}e^{-\|\bx(\bw)\|^2}.
\end{equation}
\end{lemma}
In this formula and below, we identify $\Tor$ with segment $[0,1)$ equipped with addition modulo one, and thus, the characters $\chi(\bw)$ are (after this identification) real valued.

{In particular, the proof of Lemma \ref{LmThetaLike} shows that for almost every $\cL$, whenever $\bw \neq 0$, $y(\bw) \neq 0$ (see Section \ref{ScASConv}) and that each individual summand in \eqref{Lim2k}  is finite almost everywhere on $\Phase.$}  In order to simplify the notation, we will abbreviate expressions such as \eqref{Lim2k} by
\begin{equation}
\label{DefCL-M}
\cX(\cL,\chi)=
\sum_{\bw\in\cL \setminus \{\boldsymbol{0}\}}
\frac{\sin(2\pi\chi(\bw))}{y(\bw)}e^{-\|\bx(\bw)\|^2}
\end{equation}
even though \eqref{DefCL-M} does not converge absolutely.

If we assume that the pair $(\cL, \chi)$ is distributed according to the Haar measure on $\Phase$
 then $\cX$, { defined by \eqref{Lim2k}}, becomes a random
variable. This is the variable mentioned in Theorem \ref{Th2EdgeRand}. Note that the distribution of $\cX$ depends neither on $\bP$ nor on~$z.$

{Next, we describe how we can use the second representation of Haar measure to describe $\cX$.} 
Let $\bw_1,\dots, \bw_d$ be the shortest spanning set of $\mathcal{L}$, i.e., $\bw_1$ is the shortest non zero vector in $\cL$ and, for $j>1,$ $\bw_j$ is the shortest vector in $\cL$ that is linearly
independent of $\bw_1,\dots, \bw_{j-1}.$
Given $\bm=(m_1, \dots, m_d)\in \integers^d$, {let
\begin{equation}\label{MthPoint}
  (y,\bx)(\bm):=m_1 \bw_1+\dots+m_d \bw_d\in \mathcal{L}
\end{equation}
where
$y\in \reals$ and $\bx\in \reals^{d-1}.$}
Let $\theta_j=\chi(\bw_j).$ Then $\theta_j$ are uniformly distributed on $\Tor$ and independent of each other. Set {
\begin{equation}\label{eqTheta}
  \theta(\bm):=m_1\theta_1+\dots+m_d\theta_d.
\end{equation}}
Now, $\cX$ (see definition in Lemma \ref{LmThetaLike}) can be rewritten as
\begin{equation}
\label{cXTheta}
\cX=
\sum_{\bm\in\integers^d \setminus \{\boldsymbol{0}\}}
\frac{\sin(2\pi\theta(\bm))}{y(\bm)}e^{-||\bx(\bm)||^2}
\end{equation}
where $\cL$ is uniformly distributed on the space of lattices,
$(y, \bx)(\bm)$ is defined by \eqref{MthPoint},
 and $(\theta_1, \dots \theta_d)$ is uniformly distributed on
$\Tor^d$ and independent of $\cL.$ {We will use the representation \eqref{cXTheta}
in Sections \ref{HomFlow} and \ref{ScSegm} in our proofs and in Section \ref{ScASConv} when establishing the convergence of $\cX$.}


Theorems \ref{ThDioEdge} and \ref{Th2EdgeRand} have analogues in case we are interested in probability that $S_n$ belongs to a finite interval. In particular, our results have applications to Local Limit Theorems.
\begin{theorem}
\label{ThJoint}
  Let $z_1(n)$ and $z_2(n)$ be two uniformly bounded sequences such that $|z_1(n)-z_2(n)|n^{d/2}\to\infty.$
  Then the random vector\vspace{3pt}
\begin{equation}
\label{2Point}
\frac{n^{d/2}}{\Lambda(\ba, \bp)}
  \left(e^{z_1^2/2} \left[\cE_d(z_1)-\Prob_{\ba, \bp} \left(\frac{S_n}{\sigma \sqrt{n}}\leq z_1\right)
  \right],
  e^{z_2^2/2} \left[\cE_d(z_2)-\Prob_{\ba, \bp} \left(\frac{S_n}{\sigma \sqrt{n}}\leq z_2\right)\right]\right)
\end{equation}
  converges in law to a random vector $(\cX(\cL, \chi_1), \cX(\cL, \chi_2))$ where
$\cX(\cL, \chi)$ is defined by \eqref{DefCL-M} and the triple $(\cL, \chi_1, \chi_2)$
  is uniformly distributed on $\left(SL_d(\reals)/SL_d(\integers)\right)\times \Tor^d\times \Tor^d.$
\end{theorem}
Here and below the uniform distribution of $(\cL, \chi_1, \chi_2)$ means that $\cL$ is uniformly distributed
on the space of lattices, and for a given lattice, $\chi_1$ and $\chi_2$ are chosen independently and uniformly
from the space of characters.

\begin{theorem}
\label{ThLLTRand}
  Let $z_1(n)<z_2(n)$ be two uniformly bounded sequences such that $l_n=z_2(n)-z_1(n)\to 0.$

  $(a)$ If $l_n\geq C n^{\eps-d/2}$ for some $\eps>0$ then
$$ \frac{\Prob_{\ba, \bp}(z_1<\frac{S_n}{\sigma \sqrt n}<z_2)}{l_n \fn(z_1)}\to 1 $$
  almost surely.

  $(b)$ If   $l_n n^{d/2}\to\infty$ then
$$ \frac{\Prob_{\ba, \bp} (z_1<\frac{S_n}{\sigma \sqrt n}<z_2)}{l_n \fn(z_1)}\Rightarrow 1 $$
$($here and below ``$\Rightarrow$'' denotes the convergence in law$)$.

  $(c)$ If   $l_n =\frac{c|a_{d+1}-a_1|}{\sigma(\ba, \bp) n^{d/2}}$ then
  $\DS H(\ba, \bp)  \left[ \frac{\Prob_{\ba, \bp}(z_1<\frac{S_n}{\sigma \sqrt n}<z_2)}{l_n \fn(z_1)}-1\right]\Rightarrow \cY $
where
  $$
 {H(\ba, \bp)={(2\pi)^d} \sqrt{\det(D_{\ba, \bp})}}
$$
 and
$$\cY(\cL,\chi, c)={\frac{1}{c}
\sum_{ \bw\in \cL \setminus{\{\boldsymbol{0}\}}}
\frac{\sin(2\pi\chi(\bw))-\sin(2\pi [\chi(\bw)-c y(\bw)])}
     {y(\bw)}e^{- \|\bx(\bw)\|^2}}\,, $$
$\cL, \chi$ are as in Theorem \ref{Th2EdgeRand} and $D_{\ba,\bp}$ is from \eqref{DefLambda}.
\end{theorem}
\begin{remark}
The normalization in Theorem \ref{ThLLTRand}(c) comes from the following computation.
Denote
$\DS \Delta_n(z)=\cE_d(z)- \Prob \left(\frac{S_n}{\sigma \sqrt{n}} \leq z\right).$
Then, Theorem \ref{Th2EdgeRand} can be informally restated as
$$\Delta_n(z)\approx \frac{\Lambda(\ba, \bp) \sqrt{2\pi}\, \fn(z)}{n^{d/2}}\cX . $$
Then under the assumption of part (c) of Theorem \ref{ThLLTRand}
we have
$$\frac{\Delta_n(z_2)-\Delta_n(z_1)}{l_n}
\approx \frac{\Lambda(\ba, \bp) \sqrt{2\pi} }{l_n n^{d/2}}
[\fn(z_2)\cX_2-\fn(z_1)\cX_1] . $$
Since
$\DS \frac{\Lambda(\ba, \bp) \sqrt{2\pi}}{l_n n^{d/2}}=\frac{1}{c\, H(\ba, \bp)}$
we can rewrite the above equation as
$$ c\, H(\ba,\bp) \frac{\Delta_n(z_2)-\Delta_n(z_1)}{l_n \fn(z_1)} \approx \frac{\fn(z_2)}{\fn(z_1)}\cX_2-\cX_1. $$
Thus, the proof of Theorem \ref{ThLLTRand} 
proceeds by describing the asymptotics of  the joint distributions of $n^{d/2} \Delta_n(z_1)$ and
$n^{d/2} \Delta_n(z_2)$
while Theorem \ref{Th2EdgeRand} gives the marginal distributions.
\end{remark}

The intuition behind the results of Theorem \ref{ThLLTRand} is the following. Call $y_n$
$\delta$-plausible if $ \Prob(S_n=y_n)\geq \delta n^{-d/2}$. The discussion following \eqref{PdfSum} shows that for each $\delta$ there are about $C(\delta) n^{d/2}$ $\delta$-plausible values. Therefore, if $l_n\ll n^{-d/2}$ then the interval $[z_1(n), z_2(n)]$ would typically contain no plausible values. Hence, we should not expect a Local Limit Theorem (LLT) to hold on that scale. Theorem \ref{ThLLTRand} shows that as soon as interval $[z_1(n), z_2(n)]$ contains many plausible values then an LLT typically holds for this interval.

Recall that
$$\Prob_{\ba, \bp}(S_n \in [z_1, z_2])=
\sum_{{m_i \geq 0, \;\sum m_i=n}\atop{z_1\leq \sum m_ia_i\leq z_2}}
\frac{n!}{m_1! \dots m_{d+1}!} p_1^{m_1} \dots  p_{d+1}^{m_{d+1}}\,. $$
So, in Theorem \ref{ThLLTRand}, we just count the number of visits of a random linear form $\sum m_i a_i$
to a finite interval with weights given by multinomial coefficients. It is also interesting to consider counting with
equal weight. In this case the analogue of Theorem \ref{ThLLTRand}(c) is obtained in \cite{Mar1}
while for longer intervals only
partial results are available, see \cite{DFSurv, Ke}.

The layout of the paper is the following. Theorem \ref{ThDioEdge} is proven in Section~\ref{ScDio}. The proof is a minor modification of the arguments of \cite[Chapter XVI]{Fel}. The bulk of the paper is devoted to the proof of
Theorem~\ref{Th2EdgeRand}. In Section \ref{ScChange},
we provide an equivalent formula for $\cX.$ This formula looks more complicated
than \eqref{DefCL-M} but it is easier to identify with the limit of the error term.
Section \ref{ScSmoothing} contains preliminary reductions.
 Namely, we show that
the integration in the Fourier inversion formula could be restricted to a finite domain.
In Section \ref{SS53},
we show that the main contribution to the error term comes from resonances where the characteristic function
of $S_n$ is close to 1 in absolute value.
The proof relies on the asymptotic analysis of the resonant term performed in Section \ref{ScSimplify}.
Several technical estimates used in our analysis  are established in Section~\ref{ScExpChar}.
In Section~\ref{HomFlow}, we use dynamics on homogenuous spaces in order to show that the contribution of resonances
converges to \eqref{DefCL-M} completing the proof of Theorem \ref{Th2EdgeRand}.
The proofs of Theorems \ref{ThJoint} and \ref{ThLLTRand}  are similar to the proof of Theorem \ref{Th2EdgeRand}.
The necessary modifications are explained in Section \ref{ScSegm}. Finally, Section \ref{ScASConv} contains the proof of Lemma \ref{LmThetaLike}.

{As a notational remark, in the paper the constants denoted by $c$, $C$, or other implied constants may change from line to line or even within the same line.}

\section{Edgeworth Expansion under Diophantine conditions.}
\label{ScDio}

\Cref{ThDioEdge} is a consequence of \Cref{IneqChar} and \Cref{ThDioEdgeGen} below.

Note that the characteristic function
of $X$ is given by
\begin{equation}
\label{AtomicCF}
\phi(s)=p_1 e^{is a_1}+\dots+p_{d+1} e^{is a_{d+1}}
\end{equation}
and recall that
$\DS d(s)=\max_{j\in \{2,\dots, d+1\}} \dist(b_j s, 2\pi \integers)$
where $b_j=a_j-a_1$.
{
\begin{lemma}\label{IneqChar}
There exists a positive constant $c$ such that
\begin{equation}
\label{KeyIneqChar}
|\phi(s)|\leq 1-c\, d(s)^2.
\end{equation}
\end{lemma}
}
\begin{proof}
{ Since
$$ \DS  1-|\phi(s)|=\frac{1-|\phi(s)|^2}{1+|\phi(s)|}\geq
\frac{1-|\phi(s)|^2}{2},$$
it suffices to show that
\begin{equation}
\label{KeyIneqCharSq}
|\phi(s)|^2\leq 1-2c\, d(s)^2.
\end{equation}
Note that
$$ |\phi(s)|^2=\sum_j p_j^2+2\sum_{j<k} p_j p_k \cos((a_j-a_k)s). $$
Taking a constant $\brc$ such that $\cos(t) \leq 1-\brc\, t^2$ for $|t|\leq \pi$ and letting
$c=\brc (\min_j p_j)^2$ we obtain
$$ |\phi(s)|^2\leq 1-{2c}\sum_{j<k} \text{dist}^2((a_j-a_k)s, 2\pi \integers)$$
proving \eqref{KeyIneqCharSq}.
}
\end{proof}

\begin{theorem}
\label{ThDioEdgeGen}
If the distribution of $X$ has $d+2$ moments and its characteristic function $\phi$ satisfies
\begin{equation}
\label{CharDioGen}
 |\phi(s)|\leq 1-\frac{K}{|s|^\gamma}
\end{equation}
and $R<\frac{d}{2}$ is such that
\begin{equation}
   \label{RSmallGen}
   \left(R-\frac{1}{2}\right)\gamma<1
   \end{equation}
  then
\begin{equation}\label{Edge(d-1)}
\lim_{n\to\infty} n^R \left[{\Prob} \left(\frac{S_n}{\sigma \sqrt{n}}\leq z\right)-\cE_{d-1}(z)\right]
=0.
\end{equation}

{ In particular, if $X$ is discrete with $d+1$ atoms
$\ba=(a_1, \dots, a_{d+1})\,,$ $\ba$ is $\beta-$Diophantine
and $ \left(R-\frac{1}{2}\right)\beta<\frac{1}{2}\,, $ then \eqref{Edge(d-1)} holds.}

\end{theorem}

Theorem \ref{ThDioEdgeGen} follows easily from the estimates in \cite[ChapterXVI]{Fel} but we provide the proof here for
completeness.
\begin{proof}
Denoting
$${ \brDelta_n=\Prob}\left(\frac{S_n}{\sigma \sqrt{n}}\leq z\right)-\cE_{d-1}(z) $$
we get, {from the estimate (3.13) in \cite[Chapter XVI]{Fel}}, that for each $T$
\begin{equation}\label{Estr}
  |\brDelta_n|\leq
  \frac{1}{\pi}\int_{-\frac{T}{\sigma \sqrt{n}}}^{\frac{T}{\sigma \sqrt{n}}}\left|\frac{\phi^n(s)-\hcE_{d-1}(s\sigma\sqrt{n})}{s}\right|\,ds+\frac{C}{T}.
\end{equation}
Choose $T=Bn^{R}$ with $B=\frac{C}{\varepsilon}$.
Then, $\frac{C}{T}= \frac{\varepsilon }{n^{R}}$.
Take a small $\delta$ and split {the integral in the RHS} of \eqref{Estr} into two parts.
\begin{multline}\label{Splitr}
\frac{1}{\pi}\int_{-\delta}^{\delta}\left|\frac{\phi^n(s)-\hcE_{d-1}(s\sigma\sqrt{n})}{s}\right|\,ds +\frac{1}{\pi}\int_{\delta<|s|<Bn^{R-1/2}/\sigma} \left|\frac{\phi^n(s)-\hcE_{d-1}(s\sigma\sqrt{n})}{s}\right|\,ds.
\end{multline}
{From the proof of Theorem 2 in Section 2 and Theorem 3 in Section 4 of \cite[Chapter XVI]{Fel}}, we have that the first integral of \eqref{Splitr} is
$\cO\left(n^{-d/2}\right)$.

Also, $\int_{|s|>\delta}\left|\frac{\hcE_{d-1}(s\sigma \sqrt{n})}{s}\right| \, ds$
has exponential decay as $n \to \infty$. Put $$J=\{s:\;\delta<|s|<B\sigma^{-1}n^{R-1/2}\}.$$
Thus, we only need to estimate
\begin{equation}\label{CharInt}
  \int_{J} \left|\frac{\phi^n(s)}{s}\right| \, ds \leq
  \frac{1}{\delta}\int_{J} \left|\phi^n(s)\right| \, ds \leq
  \frac{C}{\delta}\int_{J}
  \exp\left(-b\, n^{1-\left(R-\frac{1}{2}\right)\gamma}\right)
  \, ds
\end{equation}
where the last inequality is due to \eqref{CharDioGen}.
By \eqref{RSmallGen} the integral decay faster than any power of $n.$ 
The result follows.
\end{proof}

\begin{remark}
{The fact that the Edgeworth expansion of order $(d-1)$ holds for almost every $\ba$ is obtained
in \cite[Section 4]{AP} (with a weaker error bound). \cite{B} shows, among other things, that a Diophantine condition with any exponent is sufficient for obtaining a weak Edgeworth expansion for sufficiently smooth functions. \cite{CP, FL}  obtain similar results for dependent random variable including finite state Markov chains. The relation between the Edgeworth expansions and Diophantine approximations are utilized in \cite{AP, Bo14, Bo19} to show that Edgeworth expansions hold for almost every member of several multi-parameter families.}
\end{remark}


\section{Change of variables.}
\label{ScChange}
Here, we deduce Theorem \ref{Th2EdgeRand} from:

{\bf Theorem \ref{Th2EdgeRand}*}
{\it
For each $z$ the random variable
$$ n^{d/2}
\left[\cE_d(z)-\Prob_{\ba, \bp}\left(\frac{S_n}{\sigma \sqrt{n}}\leq z\right)\right]$$
converges in law to $\hat\hcX$ where
$$ \hat\hcX(\fa,\fp,\cL,\chi)=$$
\begin{equation}
\label{DefhCX}
e^{-z^2/2}
\frac{|\fa_{d+1}-\fa_1|}{2\sigma(\fa,\fp)\sqrt{2\pi^3}}
\sum_{\bw\in \cL \setminus{\{\boldsymbol{0}\}}} \frac{\sin 2\pi\chi(\bw)}
    {y(\bw)}e^{-4\pi^2\bx(\bw)D_{\fa,\fp}\cdot \bx(\bw)}\,,
\end{equation} $(\cL, \chi)$ is distributed according to {$\mu,$ the Haar measure on $\Phase ,$}}
$\fa=(\fa_1,\dots, \fa_{d+1}),$ $\fp=(\fp_1, \dots, \fp_{d+1})$ and $(\fa,\fp)\in \Omega$
are distributed according to $\bP,$
$(\fa, \fp)$ and $(\cL, \chi)$ are independent,
and $D_{\fa, \fp}$ and $\sigma(\fa, \fp)$
are defined immediately after \eqref{DefLambda}.

{We note that the convergence of \eqref{DefhCX} for almost every
$(\cL, \chi)$ follows\footnote{{Lemma \ref{LmThetaLike} shows that the convergence holds if the
sum in \eqref{DefhCX} is understood as a limit as $R\to\infty$ of the sums restricted to the domain $\|A\bw\|\leq R$ where $A$ is the matrix given by \eqref{DefLinCV}. However, the proof of Lemma \ref{LmThetaLike} shows that this sum could also be understood as the limit of sums over domains $\|\bw\|\leq R.$}} from Lemma~\ref{LmThetaLike}, see Step 1 in the proof of Theorem  \ref{Th2EdgeRand} below.}

\begin{proof}[Proof of Theorem \ref{Th2EdgeRand} assuming Theorem \ref{Th2EdgeRand}*]
We divide the proof into three steps.\medskip

{\em Step 1.} We will show that
$e^{z^2/2} \frac{\hat\hcX}{\Lambda(\fa, \fp)}$ has the same distribution as $\cX$ (see \eqref{Lim2k}). To this end, we rewrite the sum in \eqref{DefhCX} as\vspace{3pt}
\begin{equation}
\label{MainSum}
 \frac{1}{(2\pi)^{d-1} \det(\sqrt{D_{\fa, \fp}})}
\sum_{ \bw\in \cL\setminus\{0\}} \frac{\sin(2\pi\chi(\bw))}{y(\bw)/((2\pi)^{d-1} \det(\sqrt{D_{\fa, \fp}}))} e^{-4\pi^2  \|\sqrt{D_{\fa,\fp}}\,\bx(\bw)\|^2}.
\end{equation}
Let $A$ be the linear map such that
\begin{equation}
\label{DefLinCV}
A(y, \bx)=\left(
\frac{y}{(2\pi)^{d-1}\sqrt{\det(D_{\fa, \fp})}},\;
2\pi \sqrt{D_{\fa, \fp}} \; \;  \bx
\right).
\end{equation}
Put $(\bar\brcL, \bar\chi)=A(\cL, \chi)$.
Then, using \eqref{InvProd}, \eqref{MainSum} can be rewritten as:
$$
\frac{1}{(2\pi)^{d-1} \det(\sqrt{D_{\fa, \fp}})}
\sum_{{\bar\bw\in { \bar\cL}\setminus\{0\}}} \frac{\sin(2\pi\bar\chi(\bar\bw))}{y(\bar\bw)} e^{-||\bx(\bar\bw)||^2}.
$$
Since $\det(A)=1,$ the pair $(\bar\brcL, \bar\chi)$ is distributed according to the Haar measure on $\Phase.$
Thus, using \eqref{DefLambda},
$$ \hat\hcX(\fa,\fp,\cL,\chi)=e^{-z^2/2}\Lambda(\fa,\fp)\sum_{{\bar\bw\in { \bar\cL}\setminus\{0\}}} \frac{\sin(2\pi\bar\chi(\bar\bw))}{y(\bar\bw)} e^{-||\bx(\bar\bw)||^2}$$
where $(\bar\brcL, \bar\chi)$ is distributed according to the Haar measure on $\Phase.$
So, from \eqref{Lim2k}
$e^{z^2/2} \frac{\hat\hcX}{\Lambda(\fa, \fp)}$ and $\cX$
have the same distribution.
\medskip

{\em Step 2.} Denote
{
\begin{equation}
\label{DefOmegaMKappa}
\Omega^M_\kappa=\{(\ba, \bp)\in \Omega: \forall i\;\;\kappa \leq p_i, \quad |a_i| \leq M\quad\text{\and}\quad
\forall i\neq j \;\; |a_i-a_j|\geq \kappa\},
\end{equation}
$$\Delta_n=\cE_d(z)- \Prob \left(\frac{S_n}{\sigma \sqrt{n}} \leq z\right),
\quad
\tDelta_n=e^{z^2/2} \frac{\Delta_n}{\Lambda(\ba, \bp)}.
$$
We claim that it is enough to prove Theorem \ref{Th2EdgeRand} under the assumption that
$\bP$ has smooth density supported  on $\Omega^M_\kappa$ for some $\kappa$ and $M.$
Indeed, let $p$ be the original density of $\bP.$ Let
$f:\reals\to\reals$ be a smooth compactly supported function.
Given $\eps$ there exists a smooth density $\tp$ supported on some $\Omega^M_\kappa$ such
that $\|\tp-p\|_{L^1}\leq \frac{\eps}{2\|f\|_\infty}. $

If Theorem \ref{Th2EdgeRand} holds for smooth
compactly supported densities then we can find $n_0\in \naturals$ such that for $n\geq n_0$
$$ \left|\iint f\left(n^{d/2} \tDelta_n\right) \tp\, d\ba\, d \bp-
\EXP(f(\cX)) \right|\leq \frac{\eps}{2}. $$
Since
$$
 \left|\iint f\left(n^{d/2} \tDelta_n\right) \tp\,  d\ba d \bp-
\iint f\left(n^{d/2} \tDelta_n\right) p\,  d\ba\, d \bp\right| \leq ||p-\tp||_{L^1} ||f||_{L^\infty} \leq \frac{\eps}{2}
$$
it follows that
$$ \left|\iint f\left(n^{d/2} \tDelta_n\right) p\, d\ba d \bp-
\EXP(f(\cX))\right|\leq \eps. $$
Since $\eps$ is arbitrary, Theorem \ref{Th2EdgeRand} follows for an arbitrary $L^1$ density.
\medskip

{\em Step 3.} By Step 2, we can and will assume that $(\ba, \bp)$ is distributed according to
a smooth density supported on $\Omega^M_\kappa$ for some $\kappa$ and $M.$
Let $f$ be a smooth compactly supported test function. Divide $\Omega_\kappa^M$ into
small cubes $\{Q_j\}$ such that if $(\ba_j, \bp_j)$ denotes the center of $Q_j$, then for each
$j$, each $(\ba, \bp)\in Q_j$ and each $\Delta \in \reals$ we have
$$ \left|f\left(\frac{\Delta}{\Lambda(\ba, \bp)} \right)-
f\left(\frac{\Delta}{\Lambda(\ba_j, \bp_j)} \right)\right|\leq \eps.$$
Such a partition exists since $\Lambda$ is continuous and bounded away from 0 on $\Omega^M_\kappa.$
Then
$$ \iint f\left(n^{d/2} \tDelta_n\right) p\, d\ba\, d\bp=
\iint_{\Omega_\kappa^M} f\left(e^{z^2/2} \frac{n^{d/2} \Delta_n}{\Lambda(\ba, \bp)} \right) p\, d\ba\, d\bp$$
$$=
\sum_j
\iint_{Q_j} f\left(e^{z^2/2} \frac{n^{d/2} \Delta_n}{\Lambda(\ba_j, \bp_j)} \right) p\, d\ba \, d\bp+\delta(n)
$$
where $|\delta(n)|\leq \eps$ for large $n.$

Applying Theorem \ref{Th2EdgeRand}* in the case where
$(\ba, \bp)$ is distributed according to $\bP$ conditioned on $Q_j$, we get
$$ \lim_{n\to\infty}
\iint_{Q_j} f\left(e^{z^2/2} \frac{n^{d/2} \Delta_n}{\Lambda(\ba_j, \bp_j)} \right) p\, d\ba \, d\bp
 =
\bP(Q_j) \EXP\left(f\left(e^{z^2/2} \frac{\hat \hcX}{\Lambda(\ba_j, \bp_j)}\right)\right)$$
$$ =
\iint_{Q_j} \EXP\left(f\left(e^{z^2/2} \frac{\hat \hcX}{\Lambda(\fa, \fp)}\right)\right)
p \, d\fa \, d\fp+\delta_{j} $$
where $|\delta_{j}|\leq \eps\bP(Q_j) .$

By Step 1,
$$ \iint_{Q_j} \EXP\left(f\left(e^{z^2/2} \frac{\hat \hcX}{\Lambda(\fa, \fp)}\right)\right) p\, d\fa\, d\fp=
\bP(Q_j) \EXP(f(\cX)). $$
Summing over $j$ we conclude that for large $n$
$$\left|\iint f\left(n^{d/2} \tDelta_n\right) p \, d\ba\,  d\bp-\EXP(f(\cX))\right|\leq 3\eps. $$
Since $\eps$ is arbitrary, Theorem \ref{Th2EdgeRand} follows.
}
\end{proof}

\begin{remark}
The argument of Step 3 provides the following extension of Theorem~\ref{Th2EdgeRand}:

{\it The triple $\big(\,n^{d/2} \tDelta_n(\ba, \bp), \ba, \bp\,\big)$ converges in law as $n\to\infty$
to the triple $(\cX, \fa, \fp)$ where $\cX$ has the distribution described in Theorem \ref{Th2EdgeRand},
and $(\fa, \fp)$ is distributed according to $\bP$ and is independent from $\cX.$}
\end{remark}

\begin{remark}
\label{RmDensityCutOff}
The argument of Step 2 shows that it suffices to prove Theorem \ref{Th2EdgeRand}* in the case
$\bP$ has smooth density supported on $\Omega^M_\kappa$ for some $\kappa$ and $M.$
\end{remark}

Sections \ref{ScSmoothing}--\ref{HomFlow} are devoted to the proof of Theorem \ref{Th2EdgeRand}*. Note that similarly to
\eqref{cXTheta} we have
$$ \widehat\cX=e^{-z^2/2}
\frac{|\fa_{d+1}-\fa_1|}{2\sigma(\fa,\fp)\sqrt{\pi^3}}
\sum_{\bm\in \integers^d\setminus{\{\boldsymbol{0}\}}} \frac{\sin 2\pi\theta(\bm)}
    {y(\bm)}e^{-4\pi^2\bx(\bm)D_{\fa,\fp}\cdot\bx(\bm)} . $$

\begin{remark}
    The statements of Theorems  \ref{Th2EdgeRand}   and \ref{Th2EdgeRand}* look similar,
    however, there is an important distinction.
    Namely, the proof of Theorem~\ref{Th2EdgeRand}* is constructive.
    In the course of the proof, given $n$, $\ba$ and $z$,
we  construct a lattice $\cL(\ba, n)$ and a character $\chi(\ba, \bp, n, z)$ such that the expression
$n^{-d/2} \hcX(\ba,\bp,\cL(\ba, n), \chi(\ba, \bp, n, z))$
well-approximates  the error in the Edgeworth expansion.
We believe that such a construction could be made for more general distributions where the Edgeworth expansion
fails, and this will be a subject of a future investigation. So, the difference between Theorems \ref{Th2EdgeRand}  and \ref{Th2EdgeRand}* is that in the first case, we have only an approximation in law, while in the second case, we are able
to obtain an approximation in probability.
\end{remark}

\section{Cut off.}
\label{ScSmoothing}
Here we begin the proof of Theorem \ref{Th2EdgeRand}$^*$.
By Remark \ref{RmDensityCutOff}, we may and will assume that $\bP$ has a smooth
density supported on $\Omega_\kappa^M$ for some $\kappa$ and $M.$ Moreover, all constants, including the implied ones in $\cO$-estimates, may depend on $d, \kappa, M$ and $\bP$.

As in the previous section, let
$$\Delta_n=\cE_d(z)-F_n(z)\quad\text{where}\quad
F_n(z)=\Prob_{\ba, \bp} \left(\frac{S_n}{\sigma \sqrt{n}}\leq z\right). $$
Denote by $v_T(x)=\frac{1}{\pi}\cdot \frac{1-\cos Tx}{T x^2}$
and let
$\cV(s, T)=\left(1-\frac{|s|}{T}\right) \mathbbm{1}_{|s|\leq T}$
be its Fourier transform.\footnote{{We use $\int e^{is x}f(s)\, ds $ as definition of the Fourier transform of $f\in L^1$ as in \cite[Chapter XVI]{Fel}}.}

We use the approach of \cite[Section XVI.3]{Fel}. Let $T_2=\sigma n^{2d+6} $.
Note that
$\sigma =\sigma(\ba,\bp)$ is random. Since we assume that $(\ba,\bp)\in \Omega^M_\kappa$, $\sigma$ is uniformly bounded, and bounded away from $0$. So, $T_2=\cO(n^{2d+6})$ uniformly in $(\ba,\bp)$, i.e., there exist constants $c, C>0$ such that
$\DS \lim_{n \to \infty} T_2/n^{2d+6} \in (c,C)$. 

Decompose
\begin{equation}
\label{ConvApp}
- \Delta_n= { \left[F_n-\cE_d\right]\star v_{T_2} (z)}+
\left[F_n-F_n\star v_{T_2}\right](z)-\left[\cE_d-\cE_d\star v_{T_2}\right](z).
\end{equation}

To estimate the last term, we split
\begin{align} \label{SmoothConv}
\left[\cE_d-\cE_d\star v_{T_2}\right](z)&=
\int_{|x|\leq 1} \left[\cE_d(z)-\cE_d(z-x)\right] v_{T_2}(x) dx\\
&\qquad+\int_{|x|\geq 1} \left[\cE_d(z)-\cE_d(z-x)\right] v_{T_2}(x) dx. \nonumber
\end{align}

The first integral in \eqref{SmoothConv} equals to
$$ { \int_{|x|\leq 1} \cE_d'(z) \, x\, v_{T_2}(x) dx - \int_{|x|\leq 1} \frac{\cE_d''(y(z,x))}{2} x^2 v_{T_2}(x) dx} $$
$$=-\int_{|x|\leq 1} \frac{\cE_d''(y(z,x))}{2} \left(\frac{1-\cos T_2 x}{{ \pi }T_2}\right) dx=\cO\left(\frac{1}{T_2}\right)=\cO(n^{-(2d+6)})$$
where the first equality uses that $v_T$ is even.

Since both $\cE_d$ and cosine are bounded the second integral in \eqref{SmoothConv} is
bounded by
$$ C\int_{|x|\geq 1} \frac{dx}{T_2 x^2}=\frac{C}{T_2}=\cO(n^{-(2d+6)}). $$
Thus, the last term in \eqref{ConvApp} is $\cO(n^{-(2d+6)})$. {Here and below, the constant $C$ do not depend on the choice of $(\ba,\bp)$}. 

To estimate the second term in \eqref{ConvApp}, we split the integral in
$F_n\star v_{T_2}$ into regions  $\{|x|\geq 1/\sqrt{T_2}\}$ and $\{|x| \leq 1/\sqrt{T_2}\}.$
The contribution of $\{|x|\geq 1/\sqrt{T_2} \}$ is bounded by
$$ C\int_{1/\sqrt{T_2}}^\infty \frac{dx}{T_2 x^2}=\frac{C}{\sqrt{T_2}}=\cO(n^{-(d+3)}). $$

On the other hand
$$ \int_{|x|\leq 1/\sqrt{T_2}} \left[F_n(z)-F_n(z-x)\right] v_{T_2}(x) dx=0$$
unless there is a point of increase of
$F_n$ inside the interval $$J_2=\left[z-1/\sqrt{T_2}, z+1/\sqrt{T_2}\right].$$

The probability that $J_2$ contains a point of increase of $F_n$  is bounded by\vspace{3pt}
\begin{equation}
\label{Jump-Short}
\sum_{m_1+\dots+m_{d+1}=n} \bP(\fB_{\bm})
\end{equation}
where
$$ \fB_{\bm}:= \left\{\frac{m_1 a_1+\dots+m_{d+1}a_{d+1}}{{\sigma\sqrt{n}} }\in
\left[z-1/\sqrt{T_2}, z+1/\sqrt{T_2}\right]
\right\}\,.$$
Note that $\DS \fB_\bm=\left\{\left|\bm\cdot \ba-\sigma z\sqrt{n}\right|\leq \sigma \sqrt{n/T_2}\right\}.$
Since $\sigma$ is bounded on $\Omega^M_\kappa$ there is a constant $L=L(M,\kappa)$ such that
$\DS \fB_\bm\subset \bar\fB_\bm:=\left\{\left|\bm\cdot \ba-\sigma z\sqrt{n}\right|\leq L\sqrt{n/T_2}\right\}.$
To estimate ${\bP(\bar\fB_\bm)}$ we consider the following variables on $\Omega^M_{\kappa}$:
$$\zeta=(a_1, \dots, a_{d+1}, p_1, \dots, p_{d-1}).$$
Since $\zeta$ is distributed according to the bounded density it suffices to estimate the Lesbegue measure of
$\bar\fB_{\bm}$ in these coordinates. Without loss of generality we may assume that $m_1$ is the
maximal among $(m_1, \dots, m_{d+1})$, whence $m_1>n/(d+1).$ Then for large $n$ we have
that
$$\left|\frac{\partial}{\partial a_1} \left[\bm\cdot \ba-\sigma z \sqrt{n}\right]\right|=
\left|\left[m_1-z \sqrt{n} \frac{\partial \sigma}{\partial a_1}\right]\right|\geq \frac{n}{2d}.$$
The last inequality follows because the second term is $\cO(\sqrt{n})$. Accordingly, for each fixed value of $(a_2,\dots, a_{d+1}, p_1, \dots, p_{d-1})$ the measure
of $a_1$ such that $\zeta\in \bar\fB_\bm$ belongs to the segment of length $\cO(\sqrt{n/T_2})$
is  $\cO(\sqrt{1/nT_2})$.
Hence, each term in \eqref{Jump-Short}
is $\cO\left(\frac{1}{{\sqrt{n T_2}}}\right)$, and so, the sum is
$\cO\left(\frac{n^{d}}{{\sqrt{nT_2}}}\right).$ Thus, with probability $1-\cO\left(\frac{1}{{n^{7/2}}}\right)$, we have
that $-\Delta_n=\Delta_{n,2}+{ \cO}\left(T_2^{-1/2}\right)$ where
\begin{align*}
 \Delta_{n,2}&=\frac{1}{2\pi} \int_{-T_2}^{T_2}
\frac{\left[\phi^n\left({ \frac{t}{\sigma\sqrt{n}}}\right)-\hcE_d(t)\right]}{it} \cV(t,T_2) e^{-it z} dt \\ &=
\frac{1}{2\pi} \int_{-\frac{T_2}{\sigma\sqrt{n}}}^{\frac{T_2}{\sigma \sqrt{n}}}
e^{-isz \sigma \sqrt{n}} \; \frac{\phi^n(s) -\hcE_d(s\sigma \sqrt{n})}{i s} \cV(s, n, T_2) ds \, ,
\end{align*}
$\cV(s, n, T)\!\!=\!\!1-\left|\frac{s\sigma \sqrt{n}}{T}\right|$ and $\phi(s)$ is the characteristic function of $X$ given by \eqref{AtomicCF}.

\noindent
Let $T_1=\sigma K_1n^{d/2}$ for some constant $K_1>0$, and define
$$\Delta_{n,1}
=\frac{1}{2\pi} \int_{-\frac{T_1}{\sigma \sqrt{n}}}^{\frac{T_1}{\sigma \sqrt{n}}}
e^{-isz \sigma \sqrt{n}} \; \frac{\phi^n(s) -\hcE_d(s\sigma \sqrt{n})}{i s} \; \cV(s, n, T_2)\; ds . $$
Note that $T_1 = \cO(n^{d/2})$ with the implied constant independent of $(\ba, \bp)\in \Omega_\kappa^M.$
Let $\Gamma_n=\Delta_{n,2}-\Delta_{n,1}$. Put {$J_1 = [T_1/(\sigma\sqrt{n}), T_2/(\sigma \sqrt{n})]$} and
$$ \tGamma_n
=\frac{1}{2\pi} \int_{|s|\in J_1 
}
e^{-isz \sigma \sqrt{n}} \; \frac{\phi^n(s) }{i s} \; \cV(s, n, T_2) \;ds\,.$$
Note that, due to the exponential decay of $\hcE_d$,
\begin{align*}
|\tGamma_n - \Gamma_n| \leq C \int_{|s|\in J_1}
\frac{|\hcE_d(s\sigma \sqrt{n})|}{|s|} \; ds \leq C \int_{|s|\in J_1}
\frac{e^{ - n c s^2\sigma^2}}{|s|}\, ds \leq C \, e^{- c T^2_1} \log |T_2/T_1|.
\end{align*}
Hence, there exists $\eps >0$ such that $\Gamma_n=\tGamma_n+\cO\left(e^{-\eps T_1^2}\right)$. 

Further, note that $T_1/(\sigma \sqrt{n})$ and $T_2/(\sigma \sqrt{n})$ do not depend on
$(\ba,\bp)$. Thus, $\tGamma_n$ is an integral over the union of the two intervals $J_1$ and $-J_1$
whose lengths are independent of $(\ba,\bp)$.

The main result of Section \ref{ScSmoothing} is the following.
\begin{proposition}
\label{PrSecondCutoff}
\begin{equation}
\label{SecondCutoff}
\left\Vert \tGamma_n\right\Vert_{L^2}\leq \frac{C}{\sqrt{T_1 n^{d/2}}}.
\end{equation}
\end{proposition}

The proof of Proposition \ref{PrSecondCutoff} relies on the following estimates.

\begin{lemma}
  \label{LmExpChar}
For each integer $l$  there is a constant $C=C(l)$ such that
$$ \bE\left(\left|\phi^{n-l}(s)\right|\right)\leq \frac{C}{n^{d/2}} $$
for all $|s|\geq 1$.
\end{lemma}

\begin{lemma}
\label{LmIntNearRes-1}
If $I$ is a finite interval with length of order $1$ and
$l$ be an integer
then
\begin{equation*}
\int_{I}
\; |\phi^{n-l}(s)|  
\;ds = \cO\left(\frac{1}{\sqrt{n}}\right)
\end{equation*}
 $($where the implicit constant depends on $l$ and
 on the length of $I$ but not on its location$)$.
\end{lemma}

Lemmas \ref{LmExpChar} and \ref{LmIntNearRes-1} will be proven in Section \ref{ScExpChar} and Section \ref{ScSimplify}, respectively.

\begin{proof}
[Proof of Proposition \ref{PrSecondCutoff}]
Note that $\cV$ is an even function in $s$ and $\ol{\phi(s)}=\phi(-s)$. Therefore, the complex conjugate of $\tGamma_n$ is
\begin{align*}
\ol{\tGamma}_n
&=\frac{1}{2\pi} \int_{|s|\in 
J_1}
e^{isz \sigma \sqrt{n}} \; \frac{\ol{\phi}^n(s) }{-i s} \; \cV(s, n, T_2) \;ds \\
&= \frac{1}{2\pi} \int_{|-s|\in 
J_1}
e^{i(-s)z \sigma \sqrt{n}} \; \frac{\ol{\phi}^n(-s) }{i s} \; \cV(-s, n, T_2) \;ds \\
&= \frac{1}{2\pi} \int_{|s|\in 
J_1}
e^{-isz \sigma \sqrt{n}} \; \frac{\phi^n(s) }{i s} \; \cV(s, n, T_2) \;ds = {\tGamma}_n.
\end{align*}
To estimate the $L^2$-norm of $\tGamma_n,$ we write
$$
\bE(\tGamma_n^2)
=\frac{1}{4\pi^2}\bE\left( \int_{|s|\in 
J_1}
e^{-isz\sigma  \sqrt{n}} \; \frac{\phi^n(s) }{i s} \; \cV(s, n, T_2) \;ds\right)^2
$$
$$
=-\frac{1}{4\pi^2} \iint_{|s_1|, |s_2|\in J_1} \bE\left(e^{-i (s_1+s_2)z\sigma\sqrt{n}} \phi^n(s_1) \phi^n(s_2)   \right) {\frac{ \cV_n(s_1)\, ds_1}{s_1} \frac{\cV_n(s_2)\, ds_2}{s_2}}
$$
{where
\begin{equation}\label{eq:FourierV}
  \cV_n(s)=\cV(s,n,T_2)=1-\left|\frac{s\sigma\sqrt{n}}{T_2}\right| = \left(1-\frac{|s|}{n^{2d+\frac{11}{2}}}\right)
\end{equation}
is independent of $\sigma\,,$ and $0 \leq \cV_n \leq 1$ on $J_1$.}

We split this integral into two parts.

(1) In the region where $|s_1+s_2|\leq 1$, we use Lemma
\ref{LmIntNearRes-1}  to
estimate the integral by

\begin{align}
\label{S1S2Close}
\bE\bigg(\int_{ |s_1|\in J_1} \left|\phi^n(s_1)\right| &\int_{-1-s_1}^{1-s_1}\left|\phi^n(s_2)\right|
\frac{ds_2}{{|s_2|}} \, \frac{ds_1}{{|s_1|}} \bigg) \nonumber \\ &= \cO \left(\int_{{ |s_1|}\in J_1} \frac{1}{\sqrt{n}s^2_1}
\bE\left(\left|\phi^n(s_1)\right|\right) \, ds_1 \right).
\end{align}
Plugging the estimate of Lemma \ref{LmExpChar} into \eqref{S1S2Close} and integrating  we see
that
the contribution of the first region to $\bE(\tGamma_n^2)$ is $\cO\left(\frac{1}{T_1n^{d/2}} \right).$

(2) Consider now the region where $|s_1+s_2|\geq 1.$

Recall that on $\Omega$,
\begin{equation}
\label{DefOmega}
  p_1+\dots+p_{d+1}=1, \quad \text{and}\quad p_1 a_1+\dots+p_{d+1} a_{d+1}=0.
\end{equation}
We use the $2d$-dimensional coordinates $(a_1, \bnu)$ where
$\bnu:=(p_1,p_3,\dots,p_d, b_2,\dots, b_{d+1})$.

Then there exists a compactly supported density $\rho=\rho(a_1, \bnu)$ such that the contribution of the second region is
$$ \iint_{{|s_1|, |s_2| \in J_1}\atop{|s_1+s_2|\geq 1}}\left(\iint {g(s_1,s_2,a_1,\boldsymbol{\nu})}\,\rho \, da_1  \, d\bnu\right) {{\frac{\cV_n(s_1)\, ds_1}{s_1} \frac{\cV_n(s_2)\, ds_2}{s_2}}}$$
where
$${g(s_1,s_2,a_1,\boldsymbol{\nu})= e^{-i (s_1+s_2) z \sigma \sqrt{n}} e^{in(s_1+s_2)a_1} \psi^n(s_1) \psi^n(s_2) } .$$

To estimate this integral, we integrate by parts with respect to $a_1.$ Note that
for each $k$ we have
$$\DS e^{isn a_1} =\left[\frac{1}{isn} \frac{d}{da_1} \right]^k e^{isn a_1}.$$
Fix a large $k$ (for example, we can take $k = 8d+25$). The integration by parts amounts to applying
$\left(\frac{d}{d{a_1}}\right)^k$ to $ e^{-i (s_1+s_2) z \sigma\sqrt{n}} \rho [\psi(s_1)\psi(s_2)]^n $ which leads to terms formed by products of
$$ \left\{\left(\frac{d}{d{a_1}}\right)^{k_1}\left[ e^{-i (s_1+s_2) z \sigma \sqrt{n}}\right]\right\},\,\,\,
\left\{\left(\frac{d}{d{a_1}}\right)^{k_2}[\rho ]\right\},\,\,\,\text{and}\,\,\,\left\{\left(\frac{d}{d{a_1}}\right)^{k_3} \left[\psi(s_1) \psi(s_2)\right]^n\right\} $$
where $k_1+k_2+k_3= k.$
{Note that all of the above expressions depend implicitly on $a_1$ because
$p_2$ and $p_{d+1}$ depend on $a_1$ due to the second equation in \eqref{DefOmega}.
Rewriting that equation in the form
$$\DS a_1+\sum_{j=2}^{d+1} p_{j}b_{j}=0\ ,$$ we obtain
$\DS \frac{\partial p_j}{\partial a_1}=-1/b_j,\, j=2\;\; \mathrm{or}\;\; d+1.$
We also observe that when we integrate by parts,}
the boundary terms vanish because $\rho$ is smooth and has compact support.


{
Thus, the contribution of the above term to the integral is bounded by the expectation of
$$ C \iint_{{|s_1|, |s_2|\in J_1}\atop{|s_1+s_2|\geq 1}} \frac{ n^{(k_1/2)+k_3}}{|s_1+s_2|^{k-k_1} n^k}
{\left|\phi^{n-k_3}(s_1)\right|\; \left|\phi^{n-k_3}(s_2)\right|}
\frac{ds_1}{|s_1|} \frac{ds_2}{|s_2|} .
$$
To estimate the above integral we consider two cases, $k_1 \geq k-3$ and $k_1 < k-3$.

In the first case, we use
trivial bounds $|s_1|\geq 1$, $|s_2| \geq  1, $
$|s_1+s_2|^{k-k_1} \geq 1$ and
${\left|\phi^{n-k_3}(s_2)\right|}\leq 1,$
and Lemma \ref{LmExpChar} to estimate\footnote
{{Here we use the fact that Lemma \ref{LmExpChar}
applies to any absolutely continuous distribution of $(\ba, \bp)$. In particular, it applies to
the integration with respect to the (normalized) Lebesgue measure.}}
${\EXP(\left|\phi^{n-k_3}(s_1)\right|)}$
to obtain the upper bound:
 \begin{align*}
\frac{C}{n^{d/2+k-k_1/2-k_3}} \iint_{{|s_1|, |s_2|\in J_1}}\, ds_1\, ds_2 \leq  \frac{C|J_1|^2}{n^{(d+k-3)/2}} \leq  \frac{CT^2_2}{n^{(d+k- 1)/2}} = \frac{C T^2_2}{n^{9d/2+12}} \leq \frac{C}{\sqrt{T_1 n^{d/2}}}.
\end{align*}
Since $T_1=\cO(n^{d/2})$, $T_2=\cO(n^{2d+6})$ and $k=8d+25$, we have the last inequality.
In the second case, we observe that $|s_1+s_2|^{k-k_1} \geq {|s_1+s_2|}^3$.
We divide the integration region into two parts.

(a) $|s_1+s_2|\geq 0.1 |s_2|.$ In this case the integrand is bounded by
$$\DS \frac{C}{|s_1||s_2|^4} \left|{\phi^{n-k_3} (s_1)}\right|. $$
Using Lemma \ref{LmExpChar} to estimate the expectation of the last term and then performing the
integration, we obtain the bound
$$ \frac{n^{3/2} \ln n }{n^{d/2} T_1^3}=\frac{1}{n^{d/2} T_1} \times \frac{n^{3/2} \ln n}{T_1^2}. $$
The second factor is smaller than $1$ since $T_1^2=K_1^2\sigma^2 n^d$ and $d\geq 2$. 

(b) $|s_1+s_2|\leq  0.1 |s_2|.$ In this case $s_i$'s are of the same order:
$$ \frac{1}{2}\leq  \left| \frac{s_1}{s_2}  \right|\leq 2.$$
Accordingly, the integrand is bounded by
$$ \frac{1}{s_1^2}{|s_2+s_1|^{-3}} \left|{ \phi^{n-k_3} (s_1)}\right|\;
\left| {\phi^{n-k_3} (s_2)}\right|. $$
To perform the integration over $s_2$, we divide the domain of integration
into segments $I_l(s_1)$ of length of order 1, so that {there exists $c,C>0$ such that on $I_l$,
$$c|l|\leq |s_2+s_1| \leq C|l|.$$}
Using Lemma \ref{LmIntNearRes-1} on each segment,
we obtain
\begin{equation}
\label{S2Int}
{\int_{s_2\in J_2, |s_1+s_2|<0.1 s_2}
\frac{{|\phi^{n-k_3} (s_2)|} }{|s_1+s_2|^3} ds_2\leq }
 \sum_{l} \frac{C}{l^3 \sqrt{n}}\leq \frac{C}{\sqrt{n}}
\end{equation}
{where the constant $C$ does not depend on $s_1$.}
We now perform the integration over $s_1$. Using Lemma \ref{LmExpChar}  we bound the
expectation of the integral by
\begin{equation}
\label{S1Int}
\frac{C}{n^{d/2}} \int_{|s_1|\geq T_1/(\sigma \sqrt{n} )} \frac{ds_1}{s_1^2}
=\frac{C \sqrt{n} }{n^{d/2} T_1}.
\end{equation}
Multiplying the bounds of \eqref{S2Int} and \eqref{S1Int}, we obtain that
the integral over region (b)
is also within the bounds of Proposition \ref{PrSecondCutoff}.
}
\end{proof}

Proposition \ref{PrSecondCutoff} shows that {(by taking $K_1$ sufficiently large)}
 the contribution from $\tilde{\Gamma}_n$ to the  $L^2-$limit of $n^{d/2}\Delta_n$  can be made arbitrarily small. 
 On $ |s|\leq T_1/\sigma\sqrt{n}$, due to \eqref{eq:FourierV},
we have
$$\mathcal{V}(s,n,T_2) 
=\left(1-\frac{|s|}{n^{2d+\frac{11}{2}}}\right). 
$$ 
Hence,
${\Delta_{n,1}=\hDelta_n+o(n^{-3d/2})} 
$
where
\begin{equation}
\label{DefHDelta}
\hDelta_n:=\frac{1}{2\pi} \int_{|s|\leq T_1/\sigma\sqrt{n}}
\frac{\phi^n(s)-\hcE_d(s\sigma\sqrt{n})}{is} e^{-is z\sigma\sqrt{n}} ds.
\end{equation}
{In summary, the analysis of Section \ref{ScSmoothing} shows that
$\DS n^{d/2} \|\hDelta_n-\Delta_n\|_{L^2}\to 0$ as $n\to\infty.$ }
Hence, we only need to analyze $n^{d/2}\hDelta_n$ for large $n$.

\section{Contribution of resonant intervals.}
\label{ScSimplify}

\subsection{{Definition of resonant intervals.}}\label{ResInt}
Denote $$ s_k=\frac{2\pi k}{|b_{d+1}|}$$
and let $I_k$ be the segment of length
$\frac{2\pi }{|b_{d+1}|}$ centered at $s_k.$
Let $K_2$ be a constant such that $K_2\gg K_1.$
Due to the results of the previous section, it is sufficient to study
$$\hDelta_n=\sum_{|k|\leq {K_2n^{(d-1)/2}}} \tcI_k $$
where
$$ \tcI_k=\frac{1}{2\pi i}\int_{I_k}
e^{-isz \sigma \sqrt{n}} \; \frac{\phi^n(s) -\hcE_d(s\sigma \sqrt{n})}{s}  {\mathbbm{1}_{|s|\leq T_1/\sigma\sqrt{n}}} ds. $$
${\tcI_0=o(n^{-d/2})}
$ due to \cite[Section XVI.2]{Fel}. Next, $\hcE_d(s\sigma \sqrt{n})$ decays exponentially
with respect to $n$ outside of $I_0$. So, its contribution to $\tcI_k$ is negligible for $k\neq 0.$ Accordingly,
$$ \hDelta_n=\sum_{0<|k|\leq {K_2 n^{(d-1)/2}}} \cI_k+o\left(\frac{1}{n^{d/2}}\right)$$
where
$$\cI_k=\frac{1}{2\pi i} \int_{I_k}
e^{-is z\sigma \sqrt{n}} \; \frac{\phi^n(s)}{s}  \mathbbm{1}_{|s|\leq T_1/\sigma\sqrt{n}} \;ds. $$
Write
$$
\brs_k=\arg\max_{s\in I_k} |\phi(s)|, \quad \phi(\brs_k)=r_k e^{i\phi_k}.
$$

{Call the interval $I_k$ {\it resonant} if $r_k^n\geq n^{-100d}$ and call it {\it {non-resonant}} otherwise.
By definition, if the interval $I_k$ is non-resonant, then $\cI_k=\cO(n^{-100 d})$. Since there are $\cO(n^{(d-1)/2})$ number of intervals (both resonant and non-resonant), the total contribution of the non-resonant intervals is at most $\cO(n^{-(199d+1)/2})$ which is negligible.
So, from now on, we focus only on the contribution of the resonant intervals.}

\subsection{{Asymptotics of the resonant terms.}}\label{AsymRes}
The following lemma is similar to the results of \cite[Section 5.2]{D}.
\begin{lemma}
\label{LmIntNearRes}
Suppose that
\begin{equation}
\label{BigRk}
r_k^n\geq n^{-100{d} }
\end{equation}
and
\begin{equation}
\label{NotBoundary}
\pm \frac{T_1}{\sigma\sqrt{n}}\not\in I_k.
\end{equation}
Then
$$\cI_k=\frac{1}{i\sqrt{2\pi n} \sigma} \frac{r_k^n}{\brs_k} e^{-z^2/2} \;
e^{i n\phi_k-i \brs_k z\sigma\sqrt{n}}\left(1+\cO\left({\frac{\ln^3 n}{\sqrt{n}}}\right)\right).
 $$
\end{lemma}
\begin{proof}
{
Let $e^{i\brs_k a_j}=e^{i \left(\phi_k+\beta_j(k)\right)}$ with $|\beta_j(k)| \leq \pi$. Then
\begin{equation}
\label{R-Beta}
r_k=\sum_{j=1}^{d+1} p_j \cos\beta_j(k),
\end{equation}
and
\begin{equation}
\label{SumSine}
 \sum_{j=1}^{d+1} p_j \sin \beta_j(k)=0.
\end{equation}
From \eqref{BigRk}, we have
\begin{equation}\label{LowerBdrk}
r_k\geq 1-\frac{C\ln n}{n}.
\end{equation}

\eqref{LowerBdrk} along with \eqref{R-Beta} give
\begin{equation}
\label{SumBeta2}
\sum_{j=1}^{d+1} p_j \beta_j(k)^2 \leq \frac{C\ln n}{n},
\end{equation}
and hence,
$|\beta_j(k)|\leq C\sqrt{\frac{\ln n}{n}}.$
Combining this with \eqref{SumSine} we obtain
\begin{equation}
\label{SumBeta}
\sum_{j=1}^{d+1} p_j \beta_j(k)=\cO\left(\frac{\ln^{3/2} n}{n^{3/2}}\right).
\end{equation}
Next, by the definition of $\brs_k,$ $\frac{\partial}{\partial \delta}\Big|_{\delta=0}\;  \phi(\brs_k+\delta)$
is perpendicular to $\phi(\brs_k)$ and whence
\begin{equation}
\label{SumBetaA}
\sum_j p_j a_j \sin\beta_j(k)=0.
\end{equation}

Let $s \in I_k$, then $s=\bar s_k+\delta$ for some $\delta$. Using Taylor expansion,
\begin{align*}
e^{i(\brs_k+\delta)a_j}&=e^{i\phi_k}e^{i\beta_j(k)}\left(1+ia_j\delta-\frac{a_j^2 \delta^2}{2}\right) +\cO\left(\delta^3\right)\\ &=e^{i\phi_k}\bigg(\cos\beta_j(k)+i\sin\beta_j(k) + i\delta a_j\cos\beta_j(k) - \delta a_j\sin \beta_j(k)\bigg)\\ &\phantom{aaaaaaaaaaaaaaaaaaaaa}-e^{i\phi_k}\left(\cos \beta_j(k)+i\sin \beta_j(k)\right)\frac{a_j^2\delta^2}{2}+\cO\left(\delta^3\right).
\end{align*}
Thus,
$$
 \phi(\brs_k+\delta) = \sum_{j=1}^{d+1}p_je^{i(\bar s_k +\delta)a_j}=
 $$$$
 e^{i\phi_k} r_k+e^{i\phi_k}\sum_{j=1}^{d+1} p_j \cos \beta_j(k)\left(i a_j \delta -\frac{a_j^2\delta^2}{2}\right)+\cO\left(\frac{\ln^{3/2} n}{n^{3/2}}+\delta^3\right)
 $$
 \begin{equation}
 \label{RkExpansion}
 =r_k e^{i\phi_k} \left(1-\frac{\sigma^2 \delta^2}{2 } \right) +i \delta e^{i\phi_k}\sum_{j=1}^{d+1} p_j a_j  (\cos\beta_j(k)-1)
\end{equation}
 $$
 \quad - \frac{\delta^2}{2} e^{i\phi_k}\sum_{j=1}^{d+1} p_j a^2_j  (\cos\beta_j(k)-r_k)  +\cO\left(\frac{\ln^{3/2} n}{n^{3/2}}+\delta^3\right)
$$
where we have used \eqref{SumSine},
 \eqref{SumBeta}, \eqref{SumBetaA}
 as well as
$$ p_1 a_1+\dots+p_{d+1} a_{d+1}=0\quad \text{and}\quad
p_1 a_1^2+\dots+p_{d+1} a_{d+1}^2=\sigma^2. $$
The main term in \eqref{RkExpansion} is the first one since
$$|\cos \beta_j(k)-r_k | \leq |\cos \beta_j(k) - 1| + |1-r_k| = \cO\left(\beta_j(k)^2+\frac{\ln n }{n}\right).$$
Hence,  using \eqref{SumBeta2}, we obtain
\begin{align*}
\phi(\brs_k+\delta)&=r_k e^{i\phi_k} \left(1-\frac{\sigma^2 \delta^2}{2 } \right) + \cO\left((\delta+\delta^2)\frac{\ln n}{n}\right)  +\cO\left(\frac{\ln^{3/2} n}{n^{3/2}}+\delta^3\right)\\ &=r_k e^{i\phi_k} \left(1-\frac{\sigma^2 \delta^2}{2 } \right) +\cO\left(\frac{\ln^{3/2} n}{n^{3/2}}+\delta^3\right) .
\end{align*}
{In summary,}
\begin{equation}
\label{CharPert}
\phi(\brs_k+\delta)=r_k e^{i\phi_k} \left(1-\frac{\sigma^2 \delta^2}{2 } \right)
+\cO\left(\frac{\ln^{3/2} n}{n^{3/2}}+\delta^3\right).
\end{equation}

Next, split $I_k=I_k'\cup I_k''$ where $I_k'$ is the part of $I_k$ where
$\DS \left\{\left|\delta\right|\leq \frac{C\ln n}{\sqrt{n}}\right\}$ and $I_k''=I_k\setminus I_k'.$
Note that, if \eqref{LowerBdrk} holds, then Lemma
\ref{IneqChar} shows that
$\brs_k$ is close to $s_k$. So, the set
$\DS \left\{\left|\delta\right|\leq \frac{C \ln n}{\sqrt{n}}\right\}$
is completely contained in $I_k.$ Lemma \ref{IneqChar} also shows that for
$s\in I_k'',$ $|\phi(s)|^n\leq n^{-c\ln n}$. So, the contribution of $I_k''$ to $\cI_k$ is negligible.

Next,  on $I_k'$ the error term in \eqref{CharPert} is $\cO\left(\frac{\ln^3 n}{n^{3/2}} \right).$ Hence, the contribution to $\cI_k$ from $I_k'$ is
\begin{align*}
&\frac{r^n_k}{2\pi i\brs_k} e^{i(n \phi_k-\sqrt{n}\sigma z \brs_k)}
\int_{|\delta|<C\ln n/\sqrt{n}} \left(1-\frac{\sigma^2 \delta^2}{2}+\cO\left(\frac{\ln^3 n}{n^{3/2}}\right)\right)^n
\left(1+\cO(\delta)\right)
e^{-i\sigma z \delta\sqrt{n}} d\delta \\
&=
\frac{r^n_k}{2\pi i\brs_k} e^{i(n \phi_k-\sqrt{n}\sigma z \brs_k)}
\int_{|\delta|<C\ln n/\sqrt{n}} e^{-\sigma^2 \delta^2 n/2-i\sigma \delta\sqrt{n} z}e^{\cO\left(\ln^3 n/\sqrt{n} \right)}\left(1+\cO(n\delta^4 + \delta)\right)\, d\delta \\
&=\frac{r^n_k}{2\pi i\brs_k} e^{i(n \phi_k-\sqrt{n}\sigma z \brs_k)}
\left(1+{\cO\left(\frac{\ln^3 n}{\sqrt{n}}\right)}\right)\int_{|\delta|<C\ln n/\sqrt{n}} e^{-\sigma^2 \delta^2 n/2-i\sigma \delta\sqrt{n} z}\, d\delta.
\end{align*}
Making the change of variables $\sigma \delta \sqrt{n}=t$, we can rewrite the last expression as
\begin{align*}
{} &\frac{r^n_k e^{-z^2/2}}{2\pi i\brs_k \sigma \sqrt{n}} e^{i(n \phi_k-\sqrt{n}\sigma z \brs_k)} \left(1+{\cO\left(\frac{\ln^3 n}{\sqrt{n}}\right)}\right)\int_{|\delta|<C\sigma \ln n } e^{-(t+iz)^2/2}\, dt \phantom{aaaaaaaaaaaa}\\
&= \frac{r^n_k e^{-z^2/2}}{2\pi i\brs_k \sigma \sqrt{n}} e^{i(n \phi_k-\sqrt{n}\sigma z \brs_k)}
\left(1+{\cO\left(\frac{\ln^3 n}{\sqrt{n}}\right)}\right)\int_{\reals} e^{-(t+iz)^2/2}\, dt\\
&=\frac{r^n_k e^{-z^2/2}}{\sqrt{2\pi} i\brs_k \sigma \sqrt{n}} e^{i(n \phi_k-\sqrt{n}\sigma z \brs_k)}
\left(1+{\cO\left(\frac{\ln^3 n}{\sqrt{n}}\right)}\right).
\end{align*}
This completes the proof of the lemma.}
\end{proof}
\subsection{Proof of Lemma \ref{LmIntNearRes-1}} \begin{proof}
 Note that if $|\phi^{n}(s)|\leq n^{-100 d}$ then $|\phi^{n-l}(s)|\leq n^{-50 d}$ and
    if $|\phi^{n}(s)|\geq n^{-100 d}$ then $|\phi(s)|\geq 1-\frac{C\ln n}{n},$
    and  hence,
    $|\phi^{-l}(s)|\leq 2.$ Therefore
  \begin{equation}
    \label{PhiNL}
    |\phi^{n-l}(s)|\leq 2 |\phi^n(s)|+\frac{1}{n^{50 d}}.
  \end{equation}
Thus it suffices to prove the result for $l=0.$
We can cover $I$ by a finite number of intervals
$I_k$. {For the intervals where $r^n_k < n^{-100d}$, we have}
$${\sum_{r^n_k < n^{-100d}} \int_{I_k}|\phi^n(s)|\,ds}\leq C \frac{|I|}{ n^{100{d}}}.$$
For resonant intervals where $ r^n_k \geq n^{-100d}$ and $k\neq 0$, the proof of Lemma \ref{LmIntNearRes}  shows that
$$
{\sum_{r^n_k \geq n^{-100d}} \int_{I_k}|\phi^n(s)|\,ds}\leq
C \int_{|\delta|<C\ln n/\sqrt{n} }  (1-c\delta^2)^n\, d\delta +\cO\left(n^{-c\ln n}\right)
=\cO\left(\frac{1}{\sqrt{n}}\right).
$$
Finally, the case $k=0$ is analyzed in \cite[Section XVI.2]{Fel}.
\end{proof}

\section{Simplifying the error term.}\label{SS53}

As noted above, in {the resonant case}, Lemma \ref{IneqChar}
gives $d(\brs_k)\leq C\sqrt{\frac{\ln n}{n}}.$ In particular,
$\dist(b_{d+1} \brs_k, b_{d+1} s_k)\leq  C\sqrt{\frac{\ln n}{n}}$ because $b_{d+1}s_k\in 2\pi \integers$. So,
$\xi_k:=\brs_k-s_k$ satisfies $$|\xi_k|\leq C\sqrt{\frac{\ln n}{n}}.$$
Since $d(s_k)=d(\brs_k)+\cO(s_k-\brs_k)$, we also have
$d(s_k)\leq C\sqrt{\frac{\ln n}{n}}.$

Noting that $b_j s_k=\frac{2\pi k b_j}{|b_{d+1}|}$ we define
$\eta_{j,k}=\frac{2\pi k b_j}{|b_{d+1}|}+2\pi l_{j,k}$, for $j=2,\dots, d+1\,,$ where $l_{j,k}$ is the unique integer such that
\begin{equation}
\label{DefLjk}
-\pi<\frac{2\pi k b_j}{|b_{d+1}|}+2\pi l_{j,k} \leq \pi.
\end{equation}
 Then, $\eta_{d+1,k}=0$ and the foregoing discussion gives
\begin{equation}
\label{Eta-Log}
 |\eta_{j,k}|\leq C\sqrt{\frac{\ln n}{n}}.
\end{equation}
Define the random vector $$X_k=\sqrt{n} \boldsymbol{\eta}_k$$ where
$\boldsymbol{\eta}_k$ is the vector with components $(\eta_{2,k}, \dots, \eta_{d, k})$,
and let $$Y_k=\frac{k}{n^{(d-1)/2}}.$$
{Also, for the remainder of the paper, we fix a constant $\alpha$:
\begin{equation} \label{DefAlpha} \alpha=\frac{1}{2(d-1)}.\end{equation}}

The main result of Section \ref{SS53} is the following.
\begin{proposition}
\label{PrSS53Main}
Let
$\DS \tDelta_n(\delta,K):=$
$$\frac{|b_{d+1}|e^{-z^2/2}}{n^{d/2}\sigma\sqrt{2\pi^3}} \sum_{k \in S(n,\delta, K)} \frac{\sin\left(\frac{2\pi n^{d/2}}{|b_{d+1}|}(\sqrt{n}a_1-z\sigma)Y_k + (\sqrt{n}\bq + z\sigma \bomega)\cdot X_k \right)}{Y_k}\;e^{-X_kD_{\ba,\bp}\cdot X_k}$$
where
\begin{equation}
\label{S-ndk}
S(n,\delta,K)=\{\,k>0\,|\,\,\delta<Y_k<K,\,\, |Y_k|^\alpha\|X_k\|<2^{K+1}\,\}
\end{equation}
and the vectors $\bomega=(\bomega_2, \dots, \bomega_d)$ and
$\bq=(\bq_2, \dots, \bq_d)$ satisfy
\begin{equation}
\label{OmegaQ}
\bomega_m=\frac{2 \sum_{l=1}^{d+1} p_l p_m (b_l-b_m)}{\sum_{j=1}^{d+1}\sum_{l=1}^{d+1} p_l p_j (b_l-b_j)^2}\,,
\quad \bq_m=p_m\,, \quad m=2,\dots, d.
\end{equation}
Then, given $\eps$ we can find $\delta, K$ such that
$$\Prob\left(|\hDelta_n-\tDelta_n(\delta, K)|>{\eps/n^{d/2}} \right)<\eps. $$
\end{proposition}


Before proving this, we obtain an approximation for $r_k$ and use it to obtain an approximation for $\cI_k.$

\begin{sublemma}
There exists a $(d-1) \times (d-1)$ matrix $D_{\ba,\bp}$ such that
\begin{align}\label{MagMax}
r_k =
1- \boldsymbol{\eta}_k D_{\ba,\bp} \cdot \boldsymbol{\eta}_k + \mathcal{O}(\|\boldsymbol{\eta}_k\|^{3})\,.
\end{align}
\end{sublemma}
\begin{proof}
   Writing $r^2_k = \psi(\brs_k) \overline{\psi(\brs_k)}$, $\brs_k = s_k + \xi_k$ and substituting $\eta_{j,k}+b_j\xi_k$ for $b_j\brs_k$, we obtain,
\begin{multline*}
r^2_k=\sum_{j=1}^{d+1} p^2_j+2\sum_{l>j, j \neq 1}p_lp_j\cos[(b_l-b_j)\xi_k+\eta_{l,k}-\eta_{j,k}] + 2p_{d+1}p_1\cos b_{d+1}\xi_k \\ + 2\sum_{j=2}^{d}p_jp_1\cos (b_j\xi_k+\eta_{j,k}) .
\end{multline*}
Therefore,
\begin{multline*}
r^2_k=1-\sum_{l>j, j \neq 1}p_lp_j[(b_l-b_j)\xi_k+\eta_{l,k}-\eta_{j,k}]^2 - p_{d+1}p_1b^2_{d+1}\xi^2_k \\
- \sum_{j=2}^{d}p_jp_1(b_j\xi_k+\eta_{j,k})^2 +\mathcal{O}\left(\xi_k^3+\sum_{l=2}^{d}\eta_{l,k}^3\right).
\end{multline*}
Note that the implied constants here and below can be chosen to
be independent of $(\ba,\bp) \in \Omega^M_\kappa$.

Taking $\eta_{1,k}=b_1=0$, we can write the above as
\begin{multline*}
r^2_k=-\xi^2_k\sum_{l>j}p_lp_j(b_l-b_j)^2-2\xi_k\widehat{\sum} p_lp_j(b_l-b_j)(\eta_{l,k}-\eta_{j,k})\\
+1-
\widehat{\sum}
p_lp_j(\eta_{l,k}-\eta_{j,k})^2 +\mathcal{O}\left(\xi_k^3+\sum_{l=1}^{d}\eta_{l,k}^3\right)
\end{multline*}
where the sum in $\widehat{\sum}$ is taken over the pairs $(l,j)$ such that $l>j$ and $(l,j)\neq (d+1, 1).$

Since $r^2_k$ is approximated by a quadratic polynomial in $\xi_k$
(the unknown) we can approximate $\xi_k$ by determining
${\mathrm{argmax}} \; r^2_k({\xi}), $ obtaining
\begin{align}
\label{XiEta}
\xi_k&=-\frac{\widehat{\sum}\; p_lp_j(b_l-b_j)(\eta_{l,k}-\eta_{j,k})}{\sum_{l>j}p_lp_j(b_l-b_j)^2}
+\cO\left(\Vert\boldsymbol{\eta}_k\Vert^2\right)\nonumber\\
&=-\frac{\sum_{j=1}^{d+1}\sum_{l=1}^{d+1} p_lp_j(b_l-b_j)\eta_{j,k}}{ \frac{1}{2}\sum_{j=1}^{d+1}\sum_{l=1}^{d+1} p_lp_j(b_l-b_j)^2
}+\cO\left(\Vert\boldsymbol{\eta}_k\Vert^2\right).
\end{align}
We recall that $b_1=0$ and $\eta_{1,k}=\eta_{d+1,k}=0.$
Substituting back we find $r_k$ in terms of $\eta_{j,k}$ only.
{Namely,}
$$
r^2_k=1-\widehat{\sum}\; p_lp_j(\eta_{l,k}-\eta_{j,k})^2 + \frac{ \left[
\widehat{\sum} \; p_lp_j(b_l-b_j)(\eta_{l,k}-\eta_{j,k})\right]^2}{\sum_{l>j}p_lp_j(b_l-b_j)^2} + \cO\left(\sum_{l=1}^d \eta_{l,k}^3\right).
$$
Put $R=\left[\sum_{l>j}p_lp_j(b_l-b_j)^2\right]^{-1}$. Then,
$$
  r^2_k =1+\widehat{\sum} \; p_lp_j\left[p_lp_j(b_l-b_j)^2R-1\right](\eta_{l,k}-\eta_{j,k})^2
  $$$$ +
  \sum_{\substack{l>j, m>{  \brm} \\ {(l, j)\neq (m, \brm)} \\ (l,j),
  (m,{  \brm})\neq { (d+1,1)}}}
  p_lp_jp_mp_{  \brm}
  (b_l-b_j)(b_m-b_{  \brm})
  R(\eta_{l,k}-\eta_{j,k})(\eta_{m,k}-\eta_{{  \brm},k}) +  \cO\left(\sum_{l=1}^{d}\eta_{l,k}^3\right)
$$
\begin{equation}
  \label{DefD1}
  := 1-2\sum_{l,j=2}^{d} D_{l,j}(\ba,\bp)\eta_{l,k}\eta_{j,k}+ \cO\left(\sum_{l=1}^{d}\eta_{l,k}^3\right).
\end{equation}
Thus,
$$
r_k = 1-\sum_{l,j=2}^{d} D_{l,j}(\ba,\bp)\eta_{l,k}\eta_{j,k}+  \cO\left(\sum_{l=1}^{d}\eta_{l,k}^3\right)
=1- \boldsymbol{\eta}_k D_{\ba,\bp} \cdot \boldsymbol{\eta}_k + \mathcal{O}(\|\boldsymbol{\eta}_k\|^{3})
$$
where $D_{\ba,\bp}$ is the $(d-1)\times (d-1)$ matrix with
\begin{equation}
\label{DefD2}
      [D_{\ba,\bp}]_{i,j}=D_{i,j}(\ba,\bp)\,,
\end{equation}
{proving \eqref{MagMax}}.
\end{proof}

\begin{lemma}
{The matrix $D_{\ba,\bp}$ defined by \eqref{DefD2} satisfies}
 $$\mathcal{I}_k=\frac{e^{-z^2/2}}{i\sqrt{2\pi n} \sigma} \frac{(1- \boldsymbol{\eta}_k D_{\ba,\bp} \cdot \boldsymbol{\eta}_k +
        \mathcal{O}(\|\boldsymbol{\eta}_k\|^{3}))^n}{\bar{s}_k}\; e^{in \phi_k-i\overline{s}_kz\sigma\sqrt{n}}\;(1+o(1)) $$
 where $\boldsymbol{\eta}_k=(\eta_{2,k},\dots,\eta_{d,k})$.
\end{lemma}
\begin{proof}
Follows directly from Lemma \ref{LmIntNearRes} and \eqref{MagMax}. 
\end{proof}

We next consider the $\cI_k$ at the two ends. Let $\cB(\ba, \bp)$ be the contribution of these boundary terms, i.e. from $k$ such that  $\DS\pm \frac{T_1}{\sigma\sqrt{n}}\in I_k.$
By Lemma \ref{LmIntNearRes-1}, $$\DS \cB(\ba, \bp)\leq \frac{C}{T_1}.$$
Recalling that $T_1=K_1 \sigma n^{d/2}$, we see that we can make $n^{d/2}\cB(\ba, \bp)$
as small as we wish by taking $K_1$ large. So, from now on, we ignore these terms.


\begin{lemma}
\label{LmExp-Ikl}
Let
$$\cI_{k,l}=\cI_k \mathbbm{1}_{|k|^\alpha n^{1/4} \|\boldsymbol{\eta}_k\|\in [2^l, 2^{l+1}]}.$$
For all sufficiently large $K>0$, there is a constant $\tc$ such that
$$ \bE\left(\sum_{0<|k|<Kn^{(d-1)/2}} \; \; {\widehat\sum_l}  \;\; \left|\cI_{k,l}\right|\right)=\cO\left(\frac{1}{n^{d/2}}
2^{K (d-1)} \exp\left(-\tilde c2^{2K} \right) \right) $$
{where the sum in $\widehat\sum$ is over $l$ satisfying $l>K$ and $2^l < \frac{K k^\alpha \sqrt{\ln n}}{n^{1/4}}.$
}
\end{lemma}

\begin{remark}
We could restrict to $l$ satisfying $l>K$ and $2^l < \frac{K k^\alpha \sqrt{\ln n}}{n^{1/4}}$ since, by the discussion at the beginning of Section \ref{SS53},
it is enough to consider
the intervals satisfying \eqref{Eta-Log} and we can take $K> C$ where $C$ is the constant from \eqref{Eta-Log}.
\end{remark}

The proof of the above lemma will be given in \Cref{ScExpChar}.
\Cref{LmExp-Ikl}  shows that we should focus on the contribution of $\cI_{k,l}$ with
$$0<|k|<K_2n^{(d-1)/2}\quad\mathrm{and}\quad l\leq K_2.$$

Next, we prove a result that allows us to  simplify $\hDelta_n$ even further. Recall that we are dealing with resonant $k$, that is, we assume that $r^n_k\geq n^{-100d}.$
\begin{lemma}
\label{LmPhase}
$(a)$ $\overline{s}_k=s_k - \bomega \cdot \boldsymbol{\eta}_k + \mathcal{O}(\|\boldsymbol{\eta}_k\|^2)$
where $\bomega=\bomega(\ba,\bp)$ is the $1\times (d-1)$ vector
defined in \eqref{OmegaQ}.
\\
$(b)$ If $\|\boldsymbol{\eta}_k\|=\mathcal{O}\left(\frac{\ln n}{\sqrt{n}}\right)$ then $n\phi_k = ns_ka_1+np_2\eta_{2,k}+\dots+np_d\eta_{d,k}+o(1).$
\end{lemma}
\begin{proof}
Since $\overline{s}_k-s_k={ \xi_k}$ part (a) follows by \eqref{XiEta}.

{Recall that $\phi(\brs_k)=r_k e^{i\phi_k}$, and by \eqref{CharPert}
$$\phi_k=\text{arg}\ \phi(s_k)+\mathcal{O}\left(|
\bar s_k - s_k |^3+\frac{\ln^{3/2} n}{n^{3/2}}\right).$$}
Note that, $$\phi(s_k)=e^{is_ka_{1}}(p_1+p_2e^{i\eta_{2,k}}+\dots+p_{d}e^{i\eta_{d,k}}+p_{d+1}).$$
Thus,
\begin{align*}
\text{arg}(\phi(s_k)) & =s_ka_1+\tan^{-1}\left(\frac{ p_2\sin \eta_{2,k} + \dots + p_d\sin \eta_{d,k} }{p_1+p_2\cos \eta_{2,k}+\dots+p_d\cos \eta_{d,k}+p_{d+1}}\right) \\ & = s_ka_1+\sum_{l=2}^{d} p_l \eta_{l,k} + \mathcal{O}(\|\boldsymbol{\eta}_k\|^3)
\end{align*}
since the denominator in the first line is $1+\cO(\|\boldsymbol{\eta}_k\|^2).$
Part (b) now follows easily.
\end{proof}
\begin{proof}[Proof of Proposition \ref{PrSS53Main}.]
First, we show that it is enough to consider $\cI_{k,l}$ when $$\delta n^{(d-1)/2} \leq |k|<K_2n^{(d-1)/2} \quad\mathrm{and}\quad  l\leq K_2$$ for appropriately chosen $\delta$ and $K_2$.

{Recall from Section \ref{ResInt} that
$$\hDelta_n=\sum_{0<|k|\leq {K_2 n^{(d-1)/2}}} \cI_k+
o\left(\frac{1}{n^{d/2}}\right)$$
for $K_2 \gg K_1$. 
By \Cref{LmExp-Ikl}, 
the contribution of $\cI_{k,l}$ with
$$0< |k| < K_2n^{(d-1)/2}\quad\mathrm{and}\quad l > K_2$$
can be made arbitrarily small by choosing $K_2$ large.}


{Next, we claim that the distribution of $\boldsymbol{\eta}_k$ has bounded density. Since $(\ba,\bp)$ has a bounded density on\footnote{{Recall that $\Omega_\kappa^M$ is
defined by \eqref{DefOmegaMKappa}.}}
 $\Omega^M_\kappa$, the vector $$\bb=\left(\frac{b_2}{|b_{d+1}|}, \dots, \frac{b_d}{|b_{d+1}|}\right)$$ has a bounded density  on $\DS V_\kappa^M:=\left\{ (x_1,\dots,x_{d-1}) | \forall j\quad\kappa(2M)^{-1} \leq  x_j \leq {2M\kappa^{-1}}\right\}.$ Let $L$ denote the supremum of the density of $\bb.$
 Since $\boldsymbol{\eta}_k$ is obtained by rescaling $\bb$ by $2\pi k$ and taking $\hspace
{-5pt}\mod 2\pi$, its density is bounded by
\begin{equation}
\label{EtaDensity}
\frac{L}{2\pi k} \times \left\lceil \frac{4\pi M k}{\kappa} \right\rceil \leq  \frac{4M L }{\kappa}
\end{equation}
where the second factor on the LHS accounts for the multiplicity of the fractional part on $V_\kappa^M$.
Since the RHS of \eqref{EtaDensity} is independent of $k$, the claim follows.}

Next, define
$$A_1=\{(\ba,\bp)|\;\mathcal{I}_{k,l}=0 \ \forall k,l \text{ s.t. }
0<|k|< \delta n^{(d-1)/2}  \text{ and } l \leq K_2\}. $$
Then
\begin{align*}
A_1^c &=\{(\ba,\bp)|\; \exists k, l \text{ s.t. }
0<|k|< \delta n^{(d-1)/2},\, l \leq K_2,\, \mathbbm{1}_{|k|^\alpha n^{1/4} \|\boldsymbol{\eta}_k\|\in [2^l, 2^{l+1})} =1 \}\\
&=\{(\ba,\bp)|\; \exists k\text{ s.t. } 0<|k|< \delta n^{(d-1)/2} ,\; |k|^{\alpha} n^{1/4} \Vert\boldsymbol{\eta}_k \Vert < 2^{K_2+1}\}.
\end{align*}
Thus,
\begin{align}
  \label{ResNearZero}
\bP(A_1^c) &\leq \sum_{0<|k|<\delta n^{(d-1)/2}} \bP\left(|k|^\alpha n^{1/4}\|\boldsymbol{\eta}_k\| < 2^{K_2+1}\right) \\ \notag
&\leq\sum_{0<|k|<\delta n^{(d-1)/2}} \frac{C\,2^{{(K_2+1)}(d-1)}}{|k|^{(d-1)\alpha}n^{(d-1)/4}}=
\mathcal{O}\left(\sqrt{\delta}\,2^{(K_2+1)(d-1)}\right)
\end{align}
where $\alpha=[2(d-1)]^{-1}$ (see \eqref{DefAlpha}) and the probability estimate follows from $\boldsymbol{\eta}_k$ having a bounded density.

Hence, for  $K_2$ and $ \delta$
such that $\sqrt{\delta}2^{(K_2+1){(d-1)}}$ is small, we can approximate $\hDelta_n$
by the sum of $\mathcal{I}_k$'s with $\delta \leq |k|n^{-(d-1)/2} < K_2$ and
$|k|^{\alpha}n^{1/4}\|\boldsymbol{\eta}_k\| < 2^{K_2+1}$.

Combining terms corresponding to $k$ and $-k,$
we obtain the following approximation to the distribution of $\Delta_n$ for large $n$
$$\frac{|b_{d+1}|e^{-z^2/2}}{n^{d/2}\sigma\sqrt{2\pi^3}}\sum_{k \in S(n,\delta, K)}
\frac{\sin(n\phi_k-\overline{s}_kz\sigma\sqrt{n})}{Y_k}e^{-X_kD_{\ba,\bp}\cdot X_k}$$
for appropriate choices of $K$ and $\delta$, and where $S(n,\delta, K)$ is defined in \eqref{S-ndk}.
{The restriction $Y>\delta$ in $S(n, \delta, K)$ comes from \eqref{ResNearZero}, the
upper bound $Y<K$ comes from \eqref{DefHDelta}, and the restriction
$\DS |Y_k|^\alpha\|X_k\|<2^{K+1}$ comes from Lemma~\ref{LmExp-Ikl}}.
We have also used \Cref{LmPhase}(a) and the fact that $|s_k|>c\delta n^{(d-1)/2}$
in the region we consider
to replace $\bar{s}_k$ by $s_k$.

Recall (see \eqref{OmegaQ}) that $\bq:=(p_2,\dots,p_d)$.
Lemma \ref{LmPhase}(b) shows that
\begin{align*}
  n\phi_k-\overline{s}_kz\sigma\sqrt{n}&=s_k(na_1-z\sigma \sqrt{n})+n\bq \cdot \boldsymbol{\eta}_k +
  z\sigma\sqrt{n}\bomega \cdot \boldsymbol{\eta}_k+o(1)
  \\ &=\frac{2\pi n^{d/2}}{|b_{d+1}|}(\sqrt{n}a_1-z\sigma)Y_k + (\sqrt{n}\bq \KF{+} z\sigma \bomega)\cdot X_k + o(1).
\end{align*}
Therefore, for large $n$ and $K$ and
$\delta$ such that $\sqrt{\delta}2^{(K+1)(d-1)}$ is very small, the distribution of $\hDelta_n$ is well approximated by {$\tDelta_n(\delta, K)$ completing the proof of the proposition.}
\end{proof}

\section{Expectation of the characteristic function.}
\label{ScExpChar}
\begin{proof}[Proof of Lemma \ref{LmExpChar}]
{  As in the proof of Lemma \ref{LmIntNearRes-1}, the inequality \eqref{PhiNL} shows that is suffices to consider the case $l=0.$ }
Recall that $d(s)=\max_{2 \leq j \leq d+1} d(b_j s, 0)$ where the distance is computed on the torus $\reals/(2\pi \integers).$
{  Lemma \ref{IneqChar}} shows that there is a positive constant $c$ such that
\begin{equation}\label{CharBnd}
{
\left|\phi^n(s)\right|\leq e^{-c n d(s)^2}.}
\end{equation}
To prove the lemma  we decompose
$\bE\big(e^{-cn d(s)^2}\big)$
into the pieces where $d(s)\sqrt{n}$ is of order $2^l$ for some $l\leq (\log_2 n)/2$.

Since $\bP$ has a bounded density, the distribution of the
  $(b_2 s, \dots, b_{d+1} s)$ has bounded density on $\Tor^d$
  where the bound is uniform for $|s| \geq 1$.
 Hence
  $$ \bP( c_1 \leq d(s) \leq c_2 )=\cO\left(c_2^d-c_1^d\right)$$ for all $0 \leq c_1 < c_2<1$ {uniformly in $|s|\geq 1.$}
    Therefore,
\begin{align*}
\bE\left(|\phi^n(s)|\right)&\leq C \bP\left(d(s)<\frac{1}{\sqrt{n}}\right)+
C \sum_{l=0}^{(\log_2 n)/2} \bP\left(d(s)\sqrt{n}\in [2^l, 2^{l+1})\right) e^{-c 4^l} \\
&\leq \frac{C}{n^{d/2}}+C \sum_{l=0}^{(\log_2 n)/2} \frac{{  2}^{dl}}{n^{d/2}}e^{-c 4^l} \leq \frac{C}{n^{d/2}}
\end{align*}
where the constant $C$ can be chosen uniformly for all $|s|\geq 1$. This completes the proof.
\end{proof}

\begin{proof}[Proof of Lemma \ref{LmExp-Ikl}]
  Since, by \eqref{MagMax},
  $\DS r_k=1- \boldsymbol{\eta}_k D_{\ba,\bp} \cdot\boldsymbol{\eta}_k
  + \mathcal{O}(\|\boldsymbol{\eta}_k\|^{3})$
  where the implied constant is independent of $(\ba,\bp)\in \Omega^M_\kappa$
  and $\DS |k|^{\alpha}n^{1/4}\|\boldsymbol{\eta}_k\|\in [2^l,2^{l+1})$,
  we have
  $${r_k \leq  1-
  c\frac{4^{l}}{|k|^{2\alpha}\sqrt{n}}},$$
  where $c$ is independent of $(\ba,\bp)$.
  Accordingly,
  $$ r^n_k \leq C e^{-\frac{c2^{2l}\sqrt{n}}{|k|^{2\alpha}}}.$$
  Also, similarly to the proof of Lemma \ref{LmExpChar}, we get
  $$\bP(|k|^\alpha n^{1/4} \|\boldsymbol{\eta}\| \in [2^l,2^{l+1}))
  \leq \frac{C2^{l  (d-1)}}{\sqrt{|k|} n^{(d-1)/4}}.$$
Hence,
$$ \bE(|\mathcal{I}_{k,l}|) \leq \frac{Ce^{-\frac{c2^{2l}\sqrt{n}}{|k|^{2\alpha}}}}{\sqrt{n}|k|}
\frac{2^{l  (d-1)}}{\sqrt{|k|}n^{(d-1)/4}} =
\frac{C 2^{l  (d-1)} e^{-\frac{c2^{2l}\sqrt{n}}{|k|^{2\alpha}}}}{|k|^{3/2}n^{(d+1)/4}} . $$
Thus,
$$ \widehat{\sum_{l}} \bE(|\mathcal{I}_{k,l}|) \leq \frac{C2^{K  (d-1)}
e^{-\frac{c2^{2K}\sqrt{n}}{|k|^{2\alpha}}}}{|k|^{3/2}n^{(d+1)/4}}.
$$
Therefore, we need to estimate
$$\sum_{0<|k|<Kn^{(d-1)/2}} \frac{C2^{K  (d-1)}
e^{-\frac{c2^{2K}\sqrt{n}}{|k|^{2\alpha}}}}{|k|^{3/2}n^{(d+1)/4}}$$
\begin{equation}
\label{K-Sum}
=\frac{C}{n^{d/2}}\sum_{0<|k|<Kn^{(d-1)/2}} \frac{1}{|k|}
\sqrt{\frac{2^{2K  (d-1)} n^{(d-1)/2}}{|k|}}\;e^{-\frac{c2^{2K}\sqrt{n}}{|k|^{2\alpha}}}.
\end{equation}
Split the sum over
\begin{equation}
\label{ScaleK}
|k|\in \left[\frac{Kn^{(d-1)/2}}{2^{q+1}},\frac{Kn^{(d-1)/2}}{2^q}\right)
\end{equation}
for $q\in\naturals.$
Then, for a fixed $q$ we have
$$|k|^{2\alpha} = \mathcal{O}\left(\frac{K^{\frac{1}{d-1}} \sqrt{n}}{2^{\frac{q}{d-1}}}\right).$$
So, each term in the sum \eqref{K-Sum} is of order
$$\frac{2^{K{(d-1)}+(3q/2)}}{K^{3/2} n^{(d-1)/2}} \;\exp\left(-\frac{c2^{2K+\frac{q}{d-1}}}{K^{\frac{1}{d-1}}} \right).  $$
The number of terms in \eqref{ScaleK}
is ${\cO\left(\dfrac{K n^{(d-1)/2}}{2^q}\right)}$. Hence, the sum over $k$ in \eqref{ScaleK}
is
$$\mathcal{O}\left(\frac{2^{K{  (d-1)}+q/2}}{{\sqrt{K}}}\; \exp\left(-\frac{c2^{2K+\frac{q}{d-1}}}{K^{\frac{1}{d-1}}} \right)
\right)=$$
$$
\mathcal{O}\left(\frac{2^{K{  (d-1)}}}{\sqrt{K}}\; \exp\left(-\frac{c2^{2K}}{K^{\frac{1}{d-1}}} \right)
\right)\times\cO\left(2^{q/2} \exp\left[-\left(\frac{c 2^{q/2}}{\sqrt{K}}\right)^{2/(d-1)} \right]\right)
.$$
Since\footnote{ To see this, one can, for example compare the sum in \eqref{QSum} with the integral
$\DS \int_0^\infty \exp\left(-c(x/\sqrt{K})^{2/(d-1)} \right)\; dx=O\left(\sqrt{K}\right).$}
\begin{equation}\label{QSum}
\sum_q 2^{q/2} \exp\left[-\left(\frac{c 2^{q/2}}{\sqrt{K}}\right)^{2/(d-1)} \right] \leq C \sqrt{K}
\end{equation}
 we obtain the lemma upon taking ${\tc<c/ K^{\frac{1}{d-1}}}.$
\end{proof}

\section{Relation to homogeneous flows.}
\label{HomFlow}
Given $\bu\in\reals^{d-1}, v\in \reals$ consider the following function on the space $\mathcal{M}$ of
unimodular lattices in $\reals^{d}$: \vspace{3pt}
\begin{equation}\label{ApproxRV}
\cZ(L; \bu, v)=\sum_{(y,\bx)\in L\setminus{\{\boldsymbol{0}\}}} \frac{\sin 2\pi(\bu \cdot \bx+vy)}{y}
\; e^{-4\pi^2 \bx D_{\ba,\bp}\cdot \bx}
\; \mathbbm{1}_{\{\delta<y<K,\,\,2\pi\,y^{\alpha}\|\bx\|<2^{K+1}\}}.
\end{equation}
Define $ \boldsymbol{\gamma}=  \left(\frac{b_2}{|b_{d+1}|},\dots, \frac{b_d}{|b_{d+1}|}\right).$
Introduce the following matrices
$$H_{\boldsymbol{\gamma}}=\begin{pmatrix}
1 & \boldsymbol{\gamma} \\
\boldsymbol{0}^T & I_{d-1}
\end{pmatrix}, \ \ \quad \ \ G_{t}=\begin{pmatrix}
e^{-(d-1)t}  & \boldsymbol{0}\\
\boldsymbol{0}^T &  e^tI_{d-1}
\end{pmatrix} . $$
Then, we get
\begin{equation}
\label{ErrAppHom}
{n^{d/2}\tilde{\bDelta}_n(\delta,K)=\frac{|b_{d+1}|e^{-z^2/2}}{\sigma\sqrt{2\pi^3}}
\cZ({\cL(n, \ba); \bu, v})},
\end{equation}
where
\begin{equation}
\label{DefUV}
\bu=\sqrt{n}\bq + z\sigma \boldsymbol{\omega}, \quad
v=\frac{n^{d/2}}{|b_{d+1}|}(\sqrt{n}a_1 - z\sigma),
\end{equation}
$\bomega$ and $\bq$ are given by \eqref{OmegaQ},
and
$\mathcal{L}(n,\ba)$ is the unimodular lattice  $\integers^d\ H_{\boldsymbol{\gamma}}\ G_{\frac{\ln(n)}{2}}$.

{To see this, note that, for an arbitrary vector $\left(k\,, m_{2,k}\,, \dots\,, m_{d,k}\right) \in \integers^d,$
\begin{align*}
 \left(
k\,, m_{2,k}\,, \dots\,, m_{d,k}\right) \, H_{\boldsymbol{\gamma}}\ G_{\frac{\ln(n)}{2}} &= \left(
k\,, m_{2,k}\,, \dots\,, m_{d,k}\right)   \, \begin{pmatrix}
n^{-(d-1)/2} & \sqrt{n}\,\boldsymbol{\gamma} \\
\boldsymbol{0}^T & \sqrt{n}\,I_{d-1}
\end{pmatrix}\\
&=
\left(\frac{k}{n^{(d-1)/2}}, \sqrt{n}\,k\,\boldsymbol{\gamma} + \sqrt{n}\left(m_{2,k}\,, \dots\,, m_{d,k}\right)\right)\\
&=\left(Y_{k}\,,\, (2\pi)^{-1}X_k+
\sqrt{n} \left(m_{2,k}-\ell_{2,k}\,, \dots\,, m_{d,k}-\ell_{d,k}\right) \right)
\end{align*}
where $Y_k$ and $X_k$ are as in Proposition \ref{PrSS53Main}
and $l_{j,k}$ are given by \eqref{DefLjk}.}
{Note that the second term has norm at least $\sqrt{n}$ unless
$m_{j,k}=\ell_{j,k}$ for $j=2, \dots, d$.
It follows that the only term which contributes to the RHS of \eqref{ApproxRV}
is the term with $m_{j,k}=\ell_{j,k}$ justifying \eqref{ApproxRV} and \eqref{ErrAppHom}.}

Let
$\bw_j(n,\ba)=(y_j(n,\ba),\bx_j(n,\ba)),\ j=1,\dots,d$
with $y_j \in \reals$ and $\bx_j \in \reals^{d-1}$ be the shortest spanning set of $\cL(n,\ba)$. Put
$$\theta_j(n,(\ba,\bp))=\bu \cdot \bx_j(n,\ba)+v y_j(n,\ba),\ j=1,\dots,d.$$
\begin{proposition}
\label{Mix+Scale}
If $(\ba, \bp)$ is distributed according to $\bP$ then the distribution of the random vector
$$((\ba, \bp), \cL(n, \ba), \btheta(n, (\ba, \bp))) $$
converges to
$ \bP\times \mu $ as $n\to\infty,$ where $\mu$ is the Haar measure on
$$ [SL_d(\reals)/SL_d(\integers)]\times \Tor^d.$$
\end{proposition}
If we restrict our attention only to $((\ba, \bp), \cL(n, \ba))$ then the result is standard
(see \cite[Theorem 5.8]{MS1},
as well as \cite{EM, KM, Sh}).
We refer the readers to \cite[Theorem 3]{DMS}, \cite{St} and the references therein
for extensions to $ [SL_d(\reals)/SL_d(\integers)]\times (\Tor^d)^p$
under various conditions.
Our proof of Proposition \ref{Mix+Scale} follows
the approach of the proof of Proposition 5.1 in~\cite{DF}.

\begin{proof}
We need to show that for each bounded smooth test function $f$,
\begin{multline}
\label{ConvDist}
{\int_{\Omega} f((\ba,\bp),\cL(n,\ba),\boldsymbol{\theta}(n, (\ba, \bp))) \, d\bP
\to
\int_{\Omega \times \mathcal{M}\times \mathbb{T}^d} f((\ba,\bp), \cL,
\boldsymbol{\theta}) \,
d\bP \, d\cL\, d\boldsymbol{\theta}}
\end{multline}

as $n \to \infty$. Write the Fourier series expansion of $f$ {with respect to $\theta$}
\begin{equation}
\label{FourierExp}
f((\ba,\bp),\cL(n,\ba),\boldsymbol{\theta})=
\sum_{\bk=(k_1,\dots,k_d)\in \integers^d}f_{\bk}((\ba,\bp),\cL(n,\ba)) \ e^{2\pi i \bk \cdot\btheta} .
\end{equation}
Then, it is enough to prove \eqref{ConvDist} for individual terms in \eqref{FourierExp}.

If $\bk=\boldsymbol{0}$ then by \cite[Theorem 5.8]{MS1} we can conclude that
$$\int_{\Omega} f_{\boldsymbol{0}}((\ba,\bp),\cL(n,\ba)) \, d\bP
\to \int_{\Omega \times \mathcal{M}\times \mathbb{T}^d} f_{\boldsymbol{0}}((\ba,\bp), \cL) \,
d\bP \, d\cL\, d\boldsymbol{\theta} $$
as $n \to \infty$.

Next, assume that $\bk \neq \boldsymbol{0}$.
Since $\Omega$ is $2d$ dimensional, we can use
$${(a_1,\bnu):=(a_1,(p_1, p_3, \dots, p_d, b_2,\dots, b_{d+1}))}$$ as local coordinates. In these coordinates, $\cL$ is independent of $a_1$. Hence, $y_j$'s and $\bx_j$'s are independent of $a_1$.
Note that there exists a
compactly supported density $\rho= \rho(a_1,\boldsymbol{\nu})$ such that
\begin{equation}\label{termk}
J_{n,\bk} = \int f_{\bk} \ e^{2\pi i \bk \cdot \btheta} \, d\bP=
\int
\int e^{ 2\pi i \left(v\sum y_jk_j + z\sigma\sum k_j \boldsymbol{\omega} \cdot \bx_j
{+\sqrt{n}\sum k_j \bq \cdot \bx_j}
\right)} (f_{\bk} \, \rho) \,da_1 \,d\boldsymbol{\nu}
\end{equation}
where we recall that $v$ is defined in \eqref{DefUV}.
Note that
\begin{equation*}
\int_{\mathbb{T}^d\times \Omega \times \mathcal{M}} f_{\bk}
\ e^{2\pi i \bk \cdot \btheta}\  d\theta_1 \dots \, d\theta_d \, d\bP \, d\cL =0
\end{equation*}
because
$$\int_{\mathbb{T}^d} e^{2\pi i \bk \cdot \btheta }d\theta_1 \dots d\theta_d =0.$$
Therefore, it is enough to prove that $J_{n,\bk}$ converges to $0$ as $n \to \infty$.

To this end, we use integration by parts as follows. Define
\begin{align*}
g(a_1,\boldsymbol{\nu})&=e^{2\pi i\frac{ n^{(d+1)/2}\sum y_jk_j}{|b_{d+1}|}a_1} = e^{ i n^{(d+1)/2}\phi(\boldsymbol{\nu})a_1 }
\end{align*}
where $\phi(\boldsymbol{\nu})=\frac{2\pi \sum y_jk_j}{|b_{d+1}|}$, and
\begin{align*}
h(a_1,\boldsymbol{\nu}) &=(f_{\bk}\, \rho)(a_1,\boldsymbol{\nu})\, e^{-2\pi i
\left(\left(\frac{ n^{d/2}\sum y_jk_j}{|b_{d+1}|}  - 4\pi
\sum k_j \boldsymbol{\omega} \cdot \bx_j \right)z\sigma (a_1,\boldsymbol{\nu})
-\sqrt{n}\sum k_j \bq \cdot \bx_j\right)
}\,.
\end{align*}
Then, the inner integral in \eqref{termk} is  $ \int \, g(a_1,\boldsymbol{\nu})h(a_1,\boldsymbol{\nu})\, da_1 \, $.

Let $\varepsilon > 0$. On the set $Q_{\bk}=\{\phi(\boldsymbol{\nu})> \varepsilon\}$, we can write
\begin{equation*}
g(a_1,\boldsymbol{\nu}) \, da_1  = \frac{1}{i\phi(\boldsymbol{\nu})n^{(d+1)/2}} \,
d e^{ia_1n^{(d+1)/2}\phi(\boldsymbol{\nu})} .
\end{equation*}
Integrating by parts on $Q_{\bk}$
(note that $h$ has compact support)
and using trivial bounds on $Q_{\bk}^c$, we can conclude that
\begin{align*}
|J_{n,\bk}| &\leq {C\max_{\boldsymbol{\nu}} }\left|\int \frac{ e^{ia_1n^{(d+1)/2}\phi(\boldsymbol{\nu})}}{i\phi(\boldsymbol{\nu})n^{(d+1)/2}} \,
{\frac{\partial h}{\partial a_1}}(a_1,\boldsymbol{\nu}) \, da_1 \right|+ C \bP(\{\phi(\boldsymbol{\nu}) \leq \varepsilon\}) \\ &\leq \frac{{C} }{\varepsilon n^{(d+1)/2}}\int
{\max_{\boldsymbol\nu}}\left|
{\frac{\partial h}{\partial a_1}}(a_1,\boldsymbol{\nu})\right| \, d a_1+ C \bP(\{\phi(\boldsymbol{\nu}) \leq \varepsilon\})
\end{align*}
for small enough $\varepsilon$. But $${\frac{\partial h}{\partial a_1}}
(a_1,\boldsymbol{\nu})=\cO(n^{d/2}),$$ and hence, the first term is $\cO_\eps(1/\sqrt{n})$. Therefore, first taking $n \to \infty$ and then taking $\varepsilon \to 0$ we have the required result.
\end{proof}
{Recall the definitions of $(y, \bx)(\bm)$ and $\theta(\bm)$ given by \eqref{MthPoint} and \eqref{eqTheta},  respectively. With this notation,} Proposition \ref{Mix+Scale} implies that as $n\to\infty$ the distribution of $n^{d/2}\tDelta_n(\delta, K)$ converges to the distribution of\vspace{3pt}
{\begin{equation}
\label{CX-Dom}
\hat{\cX}^{(K, \delta)}(\cL, \chi):=\frac{|\fa_{d+1}-\fa_1| e^{-z^2/2}}{2\sigma(\fa,\fp)\sqrt{2\pi^3}}\sum_{\bm\in \integers^d\setminus{\{\boldsymbol{0}\}}} \frac{\sin 2\pi \theta(\bm)}{y(\bm)}e^{-4\pi^2\bx D_{\fa,\fp}\cdot \bx}\mathbbm{1}_{\cU_{K, \delta}}
\end{equation}
where\footnote{{Note that \eqref{CX-Dom} contains an additional factor of 2 in the denominator
comparing with \eqref{ErrAppHom}. This is because in \eqref{CX-Dom} the sum is over all lattice vectors
(see \eqref{DefOmegaKdelta}) while in \eqref{ErrAppHom} we only consider the vectors with
positive $y$ coordinate.}}
\begin{equation}
\label{DefOmegaKdelta}
\cU_{K, \delta}=\{\delta<|y(\bm)|<K, \; 2\pi\, |y(\bm)|^{\alpha} \|\bx(\bm)\|<2^{K+1}\}
\end{equation}
and
$(\cL, \chi)\in \Phase$ is distributed according to $\mu.$
Therefore, Theorem \ref{Th2EdgeRand}* follows from the result below.}
{\begin{lemma}
\label{LmDomains}
$\hat\cX^{(K, \delta)}$ converges in law as $K\to\infty$ and $\delta\to 0$ to the random variable
$\hat\cX$ given by \eqref{DefhCX}.
\end{lemma}}
{Lemma \ref{LmDomains}, proven in Section \ref{ScASConv}, completes the proof of
Theorem \ref{Th2EdgeRand}*.}

\section{Finite Intervals.}
\label{ScSegm}
The proofs of Theorems \ref{ThJoint} and \ref{ThLLTRand} are similar to the proofs of Theorems
\ref{ThDioEdge} and  \ref{Th2EdgeRand} so we just explain the necessary changes leaving the details to the
readers.
\begin{proof}[Proof of Theorem \ref{ThJoint}]
The random vector \eqref{2Point} can be approximated by $(\cZ^{(1)},\cZ^{(2)})$
where $\cZ^{(i)}$ are defined as in \eqref{ApproxRV} but with $\bu$ and $v$ replaced by $$\bu^{(i)} = \sqrt{n}\bq + z_{i}\sigma \boldsymbol{\omega}\ \ \text{and} \ \ v^{(i)}=\frac{n^{d/2}}{|b_{d+1}|}(\sqrt{n}a_1-z_i\sigma)$$ respectively. Define $\btheta^{(i)}$ as in Proposition \ref{Mix+Scale} but $\bu$ and $v$ replaced by $\bu^{(i)}$ and $v^{(i)}$.

To complete the proof, we prove an analogue of Proposition \ref{Mix+Scale}. Namely,
we prove that the distribution of
$$((\ba,\bp),\cL(n,\ba),\btheta^{(1)}(n, (\ba,\bp)),\btheta^{(2)}(n, (\ba,\bp)))$$ converges to $\bP\times \mu'$ as $n \to \infty$ where $\mu'$ is the Haar measure on $[SL_d(\reals)/SL_d(\integers)]\times \Tor^d \times \Tor^d$.

As in the proof of Proposition \ref{Mix+Scale}, we prove that  for individual terms in the Fourier series of a smooth function $f$ on $[SL_d(\reals)/SL_d(\integers)]\times \Tor^d \times \Tor^d$
$$\sum_{(\bk_1,\bk_2)\in \integers^d\times \integers^d}f_{\bk_1,\bk_2}((\ba,\bp),\cL(n,\ba)) \ e^{2\pi i [\bk_1\cdot \btheta^{(1)}+\bk_2\cdot (\btheta^{(1)}-\btheta^{(2)})]} $$
we have
$$ J_{n,\bk_1,\bk_2} :=
\int_{\Omega} f_{\bk_1,\bk_2}((\ba,\bp),\cL(n,\ba)) e^{2\pi i [\bk_1 \cdot \btheta^{(1)}+\bk_2 \cdot (\btheta^{(1)}-\btheta^{(2)})]} \, d\bP $$
$$ \xrightarrow{n \to \infty} \int_{\Omega \times \mathcal{M}\times \mathbb{T}^d \times \mathbb{T}^d} f_{\bk_1,\bk_2}((\ba,\bp),\cL) e^{2\pi i [\bk_1 \cdot \btheta_1+\bk_2\cdot (\btheta_1-\btheta_2)]}  \,
d\bP \, d\cL\, d\boldsymbol{\theta}_1 d\boldsymbol{\theta}_2 .
$$
The case $\bk_1=\bk_2=0$ follows from \cite[Theorem 5.8]{MS1}.
Note that
$$\bk_2 \cdot (\btheta^{(1)}-\btheta^{(2)})= (z_2(n)-z_1(n))\left(\frac{ n^{d/2}}{|b_{d+1}|}
\sum y_jk_{2,j} - \sum k_{2,j} \boldsymbol{\omega} \cdot \bx_j \right) \sigma.$$
If $\bk_1 = 0$ choose appropriate local coordinates in which $\sigma$ is a coordinate.
Integrating by parts with respect to $\sigma=\sigma(\ba,\bp)$ and using $|z_1(n)-z_2(n)|n^{d/2} \to \infty$, we see that $J_{n,\boldsymbol{0},\bk_2} \to 0$ as $n \to \infty$.

If $\bk_1 \neq 0$, then using the same local coordinates $(a_1,\bnu)$ as in the proof of Proposition \ref{Mix+Scale}, we can integrate by parts to conclude that $J_{n,\bk_1,\bk_2} \to 0$ as $n \to \infty$. The proof follows through because the leading term of $\bk_1\cdot \btheta^{(1)}+\bk_2\cdot(\btheta^{(1)}-\btheta^{(2)})$ is still $n^{(d+1)/2}\phi(\bnu)a_1$.
\end{proof}
\begin{proof}[Proof of Theorem \ref{ThLLTRand}]

 To prove part (a) pick $\bareps<\eps.$ Applying Theorem~\ref{ThDioEdge} we obtain that for almost every $(\ba, \bp)$
 $$ \Prob_{(\ba, \bp)} \left(z_1\leq \frac{S_n}{\sigma \sqrt{n}} \leq z_2\right)=
 \cE_{d-1}(z_2)-\cE_{d-1}(z_1)+\cO\left(n^{-(d-\bareps)/2}\right)$$
 $$=
 \fn(z_1) l_n+\cO(l_n^2)+\cO(l_n/\sqrt{n})+\cO\left(n^{-(d-\bareps)/2}\right) .$$
 According to the assumptions of part (a), the first term is much larger than the remaining terms proving the result.

 The proof of part (b) is similar except that we apply Theorem \ref{ThJoint} instead of
 Theorem~\ref{ThDioEdge}. So, we only get convergence in law.

 To prove part (c)  we first prove the following analogue of Theorem \ref{ThJoint} in the case where
 $z_2=z_1+\frac{c |a_{d+1}-a_1|}{n^{d/2} \sigma}$:
 $$\frac{n^{d/2}}{\Lambda(\ba, \bp)}
  \left(e^{z_1^2/2} \left[\cE_d(z_1)-\Prob_{\ba, \bp} \left(\frac{S_n}{\sigma \sqrt{n}}\leq z_1\right)\right],
  e^{z_2^2/2} \left[\cE_d(z_2)-\Prob_{\ba, \bp} \left(\frac{S_n}{\sigma \sqrt{n}}\leq z_2\right)\right]\right)$$
  converges in law to a random vector $(\tcX_1, \tcX_2)(\cL, \theta, c)$ where
$$(\tcX_1, \tcX_2)(\cL, \theta, c)=\sum_{\bm\in \integers^d\setminus{\{\boldsymbol{0}\}}}
\frac{e^{-\|\bx(\bm)\|^2}}{y(\bm)} \Big(\sin \theta(\bm), \sin(\theta(\bm)-cy(\bm))\Big). $$
Once this convergence is established, the proof of part (c) is the same as the proof of part (b).
The proof of convergence is similar to the proof of Theorem \ref{ThJoint} except that $\btheta^{(1)}$ and
$\btheta^{(2)}$ are now not independent. Namely, using the same notation as in the proof of
Theorem~\ref{ThJoint}
we have that
\begin{equation}
\label{U12V12}
{\bu^{(2)}=\bu^{(1)}+O\left(n^{-d/2}\right)\text{ and }v^{(2)}=v^{(1)}-c}.
\end{equation}
{By Proposition \ref{Mix+Scale}
$ (\cL(n, \ba),
\btheta^{(1)}(n, \ba)),$
converges as $n\to\infty$ to
$(\cL^*, \btheta^*)$
where $(\cL^*, \btheta^*)$ is distributed according to the Haar measure on
$SL_d(\reals)/SL_d(\integers)\times \Tor^d$.
Combining this fact with \eqref{U12V12} we obtain
that
$ (\cL(n, \ba),
\btheta^{(1)}(n, \ba), \btheta^{(2)}
(n, \ba)) $
converges as
$n\to\infty$
to
$(\cL^*, \btheta^*, \hat\btheta^*)$
where $(\cL^*, \btheta^*)$ is distributed according to the Haar measure on
$SL_d(\reals)/SL_d(\integers)\times \Tor^d$ and
$\hat\btheta^*_j=\btheta_j^*- c y_j.$
This justifies the formula for $(\tcX_1, \tcX_2).$}
\end{proof}

\section{Convergence of $\cX.$}
\label{ScASConv}
We need some background information. Given a piecewise smooth compactly supported
function $g: \reals^d \to \reals,$
its Siegel transform is a function on the space of lattices defined by
$$ \cS(g)(\cL)=\sum_{\bw\in \cL\setminus \{\mathbf{0}\}} g(\bw). $$
An identity of Siegel, see (\cite[Section 3.7]{Mar1} or \cite[Lecture XV]{Sie}) says that
\begin{equation}
\label{RogersEq}
 \bE_\cL(\cS(g))=\int_{\reals^d} g(\bw) d\bw.
 \end{equation}
In particular, if $B$ is a (bounded) set in $\reals^d$ with piecewise smooth boundary not containing $\mathbf{0}$ then
\begin{equation}
\label{RogersInEq}
 \bP_{\cL} (\cL\cap B\neq \emptyset)\leq \bP_{\cL}(\cS(\mathbbm{1}_B)(\cL) \geq 1)
\leq \bE_\cL (\cS(\mathbbm{1}_B))=\Vol(B).
\end{equation}

{We shall use the following consequence of this result.
\begin{lemma}
\label{LmDiopLat}
Let $\beta>d$. Then for almost every lattice $\cL$, there exist $C=C(\cL)$ such that for all
$\bw\in \cL\setminus\{\bf 0\}$ it holds
$|y(\bw)| > C\|\bw\|^{-\beta}$.
\end{lemma}
\begin{proof}
For $k\geq 1$ let $B_k\in\reals^d$ denote the following set
$$B_k=\{(\bx,y)\in \reals^{d-1}\times \reals:
\|\bx\|\in [k, k+1), |y|<k^{-\beta}\}\,, $$
and denote $D_k=\{\cL : \cL\cap B_k\neq \emptyset\}.$
By \eqref{RogersInEq}, $\bP_{\cL}(D_k)\leq K k^{d-1-\beta}.$ Since $\beta>d$,
$$\sum_k \bP_{\cL}(D_k)< \infty\,.$$ So, by the Borel-Cantelli Lemma,  for almost every
$\cL,$ there exists $k_0(\cL)$ such that $\cL\cap B_k=\emptyset$ for $k\geq k_0.$
{Since $\cL$ contains finitely many vectors satisfying $\|\bw\|< k_0$, we can choose
$C(\cL)\leq 1$
so small that $|y(\bw)|\geq C\|\bw\|^{-\beta}$ for all non-zero $\bw$ in the ball of radius $k_0.$
The result follows.}
\end{proof}
}

\begin{proof}[Proof of Lemma \ref{LmThetaLike}]
Let $\cL^+=\{\bw\in \cL: y(\bw)>0\}.$ Since
$$\frac{\sin(2\pi\chi(\bw))}{y(\bw)}$$ is even,
{and almost every lattice contains no vectors $\bw$ with $y(\bw)=0$ and $\bw \neq 0$ (this follows immediately from Lemma \ref{LmDiopLat}),}
 it is enough to restrict the attention to $\bw\in\cL^+.$ Throughout the proof we fix two numbers $\eps>0$ and $\tau<1$ such that $\eps\ll (1-\tau)\ll 1.$


By applying Lemma \ref{LmDiopLat} with $\beta=d+1$, for almost all $\cL$ we have
$$\left| \sum_{\bw\in \cL^+:\; \|x(\bw)\|\geq \|\bw\|^\eps} \frac{\sin 2\pi \chi(\bw)}{y(\bw)} e^{-\|{  x(\bw)}\|^2}\right| \leq \sum_{\bw\in \cL^+} C {\|\bw\|^{d+1}} e^{-\|\bw\|^{{2\eps}}} $$
converges absolutely. Hence, it suffices to establish the convergence of
$$ { \bar \hcX_{R}}:=
\sum_{\bw\in\cL^+: \;\|x(\bw)\|\leq \|\bw\|^\eps<R^\eps}
\frac{\sin 2\pi \chi(\bw)}{y(\bw)} e^{-\|{  x(\bw)}\|^2}. $$
Let $R_{j,k}=2^k+j 2^{\tau k},$ $j=0,\dots, \lfloor 2^{(1-\tau)k} \rfloor.$ To prove the convergence of $\hcX$,
we will show that {almost all $\cL$} and almost all $\chi$ satisfy the two
estimates below:
\begin{equation}
\label{Along}
\forall \text{ sequences } \{j_k\},\; { \bar \hcX_{R_{{j_k},k}} }\text{ converges as } k\to\infty,
\end{equation}
\begin{equation}
\label{Between}
\max_j \sup_{R_{j,k}\leq R< R_{j+1, k}} \left|\bar \hcX_{R}-{ \bar \hcX_{R_{j,k}} }\right|\to 0 \text{ as } k\to \infty.
\end{equation}
To prove \eqref{Along}, let
$$ S_{j,k}=\sum_{\bw\in\cL^+:\; \|x(\bw)\|\leq \|\bw\|^\eps, \; R_{j,k}\leq \|\bw\|< R_{j+1, k}}
\frac{\sin 2\pi \chi(\bw)}{y(\bw)} e^{-\|{  x(\bw)}\|^2}. $$
Using that
$\bE_\chi(\sin(2\pi(\chi(\bw))))=0$ and
for $\bw_1\neq \pm \bw_2$,
$$ \bE_\chi(\sin(2\pi(\chi(\bw_1)))\sin(2\pi(\chi(\bw_2))))=0,$$
we see that
$\bE_{\chi}(S_{j,k})=0$ and
\begin{align*}
\text{Var}_{\chi}(S_{j,k}) &=\sum_{\bw\in \cL^+: \; \|x(\bw)\|\leq \|\bw\|^\eps, R_{j,k}\leq \|\bw\|< R_{j+1, k}} \frac{e^{-2\|\bx(\bw)\|^2}}{2y^2(\bw)} \\
&\leq \frac{1}{2^{2k+1}} \Card(\bw: \|x(\bw)\|\leq \|\bw\|^\eps, R_{j,k}\leq \|\bw\|< R_{j+1, k}) \\
&\leq \frac{C(\cL)}{2^{2k}}\Vol(\bw: \|x(\bw)\|\leq \|\bw\|^\eps, R_{j,k}\leq \|\bw\|< R_{j+1, k})\\ &\leq \brC(\cL) 2^{(\tau+\eps(d-1)-2) k}.
\end{align*}
Hence, by Chebyshev's inequality for each $j$
$$ \bP_\chi\left(|S_{j,k}|\geq 2^{-(1-\tau+\eps)k}\right)\leq \brC(\cL) {2^{(\eps (d+1)-\tau) k}}
.$$
Therefore
$$ \bP_\chi\left(\exists j:\; |S_{j,k}|\geq 2^{-(1-\tau+\eps)k}\right)\leq \brC(\cL) {2^{(1+\eps (d+1)-2 \tau) k}} 
.$$
Thus, if $\eps$ is sufficiently small and $\tau$ is sufficiently close to $1$ then
Borel-Cantelli Lemma shows that for almost every $\chi,$ if $k$ is large enough, then for all $j$,
$|S_{j,k}|\leq 2^{-(1-\tau+\eps)k}$, and thus,
$\DS \sum_j |S_{j,k}| \leq 2^{-\eps k}$ proving \eqref{Along}.
Likewise,
$$
\sup_{R_{j,k}\leq R\leq R_{j+1, k}} \left|\bar \hcX_{R} -{   \bar \hcX_{R_{j,k}} }\right|\leq
\sum_{\bw\in\cL^+:\; \|x(\bw)\|\leq \|\bw\|^\eps, \|\bw\|\in [R_{j,k}, R_{j+1, k})}
\frac{1}{|y(\bw)|}\; e^{-\|x(\bw)\|^2}
$$$$\leq
C(\cL) 2^{-k} \Vol(\bw: \|x(\bw)\|\leq \|\bw\|^\eps, R_{j,k}\leq \|\bw\|< R_{j+1, k}) \leq
\brC(\cL) 2^{(\tau+\eps(d-1)-1)k}
$$
proving \eqref{Between}. Lemma \ref{LmThetaLike} is  established.
\end{proof}

\begin{proof}[Proof of Lemma \ref{LmDomains}]
{Given a domain $\cU\in \reals^d$ let
$$ \brcX_{\cU}(\cL, \chi)=\sum_{\bm\in \integers^d\setminus{\{\boldsymbol{0}\}}} \frac{\sin 2\pi \theta(\bm)}{y(\bm)}e^{-4\pi^2\bx D_{\fa,\fp}\cdot \bx}\mathbbm{1}_{\cU} .$$
Then
 $$\hat\cX^{(K,\delta)}=\frac{|\fa_{d+1}-\fa_1| e^{-z^2/2}}{2\sigma(\fa,\fp)\sqrt{2\pi^3}}\brcX_{\cU_{K, \delta}}$$
 where $\cU_{K, \delta}$ is given by \eqref{DefOmegaKdelta}.
Let $\Gamma_R=A B_R$ where $B_R$ is the ball of radius $R$ centered at the origin and
$A$ is the linear map given by \eqref{DefLinCV}. Lemma \ref{LmThetaLike} (after the change of variables
$\cL\mapsto A\cL$) tells us that $\brcX_{\Gamma_R}(\cL, \chi)\to \brcX(\cL, \chi)$ as $R\to\infty$ almost surely
where $\brcX=\brcX_{\reals^d}.$
Therefore it suffices to show that for each $\eta$ there exist  $\delta_0$ and $K_0$ such
if $\delta<\delta_0$ and $K\geq K_0$ then
$$ \bP\left(\left|\brcX_{\Gamma_{2^K/\delta^2}}-\brcX_{\cU_{K, \delta}}\right|>\eta\right)<\eta\,.$$
Note that $\cU_{K, \delta}\subset \Gamma_{2^K/\delta^2}$ (for sufficiently small $\delta$) and so
\begin{equation}
\label{BrCX-}
\brcX_{\Gamma_{2^K/\delta^2}}-\brcX_{\cU_{K, \delta}}=
\sum_{\bm\in \integers^d\setminus{\{\boldsymbol{0}\}}} \frac{\sin 2\pi \theta(\bm)}{y(\bm)}e^{-4\pi^2\bx D_{\fa,\fp}\cdot \bx}\mathbbm{1}_{\Gamma_{2^K/\delta^2}\setminus \cU_{K,\delta}}.
\end{equation}
Below, we choose $K$ so large and $\delta$ so small that
\begin{equation}
\label{MesFA}
\bP_{\cL}(\fA_{K, \delta}^c)\leq \eta/100
\end{equation}
 where
$\fA_{K, \delta}$ is the set of lattices $\cL$ satisfying the following conditions:}

{(i) the shortest non-zero vector in $\cL$ is longer than $3\delta^\eps$;}

{(ii) $\bw=(\bx,y)\in \cL$ then $|y|\geq \max(K, 1/\delta) \|\bw\|^{-(d+1)}$;}

{(iii) $\cL$ contains no vectors $(\bx, y)$ with $|y| \leq \delta$ and $2\pi\,\|\bx\|\leq \delta^{-1/2d}$.}

{It is easy to see using \eqref{RogersInEq} that the measure of lattices \textit{not} satisfying at least one of the above conditions
is small (cf. the proof of Lemma \ref{LmDiopLat}). We now estimate the contribution  to \eqref{BrCX-} coming from six different regions in
$\Gamma_{2^K/\delta^2}\setminus \cU_{K, \delta}.$}

{(1) Consider first the terms with $\|\bx\|\geq \|\bw\|^\eps.$ Then for $\cL\in \fA_{K, \delta}$
each term in the sum is bounded by $C \|\bw\|^{d+1} e^{-c \|\bw\|^{2\eps}}$.
We now consider several cases depending on the restrictions on $y.$}

{(a) If $|y|\geq K$ then $\|\bw\|\geq K$ and so the sum over this region is bounded
(in absolute value) by
$$
W_{1a}:=\sum_{{\bw \in \cL:\;\; \|\bw\|\geq K}} {C \|\bw\|^{d+1} e^{-c \|\bw\|^{2\eps}}}.
$$
By the Siegel identity, \eqref{RogersEq}, $\bE_\cL(W_{1a}\mathbf{1}_{\fA_{K, \delta}})\leq \brC  e^{-\brc K^{2\eps}}$
and so the contribution coming from domain (1a) is negligible in view of the Markov inequality.}


(b) If $\delta < |y|<K$, then $2\pi\,|y|^\alpha \|\bx\|\geq 2^{K+1}$ whence $\|\bx\|\geq 
2^{K}/(\pi K^\alpha).$ 
Denoting by $W_{1b}$ the contribution from this region we obtain using property (ii) of the definition
of $\fA_{K, \delta}$ that
$$
W_{1b}\leq \sum_{{\bw \in \cL:\ \|\bw\|\geq 
2^{K}/(\pi K^\alpha)}}
{C \|\bw\|^{d+1} e^{-c \|\bw\|^{2\eps}}}. $$
Hence, by \eqref{RogersEq},
$\bE_{\cL}(W_{1b}\mathbf{1}_{\fA_{K,\delta}}
) \leq
\brC e^{-\brc 2^{2\eps K}/K^{\eps\alpha}}$
which shows that the contribution from the region (1b) is  negligible in view of
the Markov inequality.

(c) If $|y| \leq \delta\,,$ then $\|\bw\|\geq \|\bx\| \geq \delta^{-1/2d}/(2\pi)$ because $\cL\in \fA_{K, \delta}.$ Hence, if $W_{1c}$ denotes
the contribution of the terms from
(1c), then,  similarly to the case (1b),
$\DS\bE_\cL(W_{1c} \mathbf{1}_{\fA_{K,\delta}})\leq \brC  e^{-\brc \delta^{-\eps/d}}$
and the contribution of region (1c) is negligible as well.

{(2) Now we discuss the terms with $\|\bx\|
\leq \|\bw\|^\eps.$ }

{Again, we shall consider three cases:}

{(a) $|y|\geq K.$
{Note that in case (2) we have $\|\bx\|^{1/\eps}\leq\|\bw\|\leq\|\bx\|+|y|$. Therefore  $|y|\geq \|\bx\|^{1/\eps}-\|\bx\| \geq (\|\bx\|/2)^{1/\eps} $ for sufficiently large $\|\bx\|$. So, $\|\bx\| \leq 2|y|^\eps.$
Thus}
taking the $L^2-$norm and integrating first with respect to $\chi$ we get
the $L^2-$norm of the terms in (2a) is bounded by
$$ \bE_\cL\left(\sum_{(\bx, y)\in \cL\cap \left(\Gamma_{2^K/\delta^2}\setminus \cU_{K, \delta}\right)
\cap\{ |y| \geq K\}}\frac{C}{y^2} \right) \leq  C\int_{\|\bx\|<2|y|^\eps, |y|\geq K}\frac{1}{y^2} d\bw  \leq \frac{C}{K^{1-\eps d}} $$
where the first inequality relies on \eqref{RogersEq}. }

{(b) $\delta<|y|<K.$ In this case,
{$\|\bx\|^{1/\eps}<\|\bx\|+K$. Thus,  $K>\|\bx\|^{1/\eps}-\|\bx\| \geq (\|\bx\|/2)^{1/\eps} $ for sufficiently large $\|\bx\|$. So,
$\|\bx\|   \leq 2K^\eps \leq K\,,$ 
and hence, $(\bx,y)\in \cU_{K, \delta}$}. So, the region (2b) does not contribute to our sum.}

{(c) $|y|\leq \delta.$ In this case, 
$\|\bw\|\leq \|\bx\|+|y| \leq  3\delta^\eps$. 
However, (ii) implies that $\|\bw\| \geq \delta^{-2/(d+1)}$. Hence, this case is impossible.
}

{Combining the six cases considered above, we obtain the result.}
 \end{proof}

\end{document}